\documentclass[12pt]{article}
\usepackage{mathrsfs}
\usepackage{amssymb}
\usepackage{amsfonts}
\usepackage{latexsym}
\usepackage{amsthm}
\usepackage{pictex}

\usepackage[cp1251]{inputenc}
\usepackage[english]{babel}

\newcommand{\er}[1]{{\rm(\ref{#1})}}
\def\lb{\label}

\theoremstyle{plain}
\newtheorem{theorem}{\bf Theorem}[section]
\newtheorem{lemma}[theorem]{\bf Lemma}

\newtheorem{proposition}[theorem]{\bf Proposition}
\theoremstyle{remark}

\begin{document}

\def\a{\alpha} 
\def\b{\beta}  
\def\g{\gamma} \def\G{\Gamma}
\def\d{\delta} \def\D{\Delta}
\def\c{\chi}
\def\z{\zeta}
\def\e{\eta}
\def\f{\phi}     \def\F{\Phi} 
\def\k{\kappa} 
\def\l{\lambda} \def\L{\Lambda}
\def\m{\mu}
\def\n{\nu}  
\def\o{\omega} \def\O{\Omega}
\def\p{\psi}   \def\P{\Psi} 
\def\r{\rho}  
\def\s{\sigma}  \def\S{\Sigma} 
\def\vT{\Theta}
\def\s{\sigma} \def\S{\Sigma} 
\def\x{\xi}   \def\X{\Xi}
\def\ve{\varepsilon}
\def\vt{\vartheta}
\def\vp{\varphi}
\def\vk{\varkappa}
\def\t{\tau}
\def\eps{\epsilon}

\def\cA{{\cal A}} \def\bA{{\bf A}}  \def\mA{{\mathscr A}}
\def\cB{{\cal B}} \def\bB{{\bf B}}  \def\mB{{\mathscr B}}
\def\cC{{\cal C}} \def\bC{{\bf C}}  \def\mC{{\mathscr C}}
\def\cD{{\cal D}} \def\bD{{\bf D}}  \def\mD{{\mathscr D}}
\def\cE{{\cal E}} \def\bE{{\bf E}}  \def\mE{{\mathscr E}}
\def\cF{{\cal F}} \def\bF{{\bf F}}  \def\mF{{\mathscr F}}
\def\cG{{\cal G}} \def\bG{{\bf G}}  \def\mG{{\mathscr G}}
\def\cH{{\cal H}} \def\bH{{\bf H}}  \def\mH{{\mathscr H}}
\def\cI{{\cal I}} \def\bI{{\bf I}}  \def\mI{{\mathscr I}}
\def\cJ{{\cal J}} \def\bJ{{\bf J}}  \def\mJ{{\mathscr J}}
\def\cK{{\cal K}} \def\bK{{\bf K}}  \def\mK{{\mathscr K}}
\def\cL{{\cal L}} \def\bL{{\bf L}}  \def\mL{{\mathscr L}}
\def\cM{{\cal M}} \def\bM{{\bf M}}  \def\mM{{\mathscr M}}
\def\cN{{\cal N}} \def\bN{{\bf N}}  \def\mN{{\mathscr N}}
\def\cO{{\cal O}} \def\bO{{\bf O}}  \def\mO{{\mathscr O}}
\def\cP{{\cal P}} \def\bP{{\bf P}}  \def\mP{{\mathscr P}}
\def\cQ{{\cal Q}} \def\bQ{{\bf Q}}  \def\mQ{{\mathscr Q}}
\def\cR{{\cal R}} \def\bR{{\bf R}}  \def\mR{{\mathscr R}}
\def\cS{{\cal S}} \def\bS{{\bf S}}  \def\mS{{\mathscr S}}
\def\cT{{\cal T}} \def\bT{{\bf T}}  \def\mT{{\mathscr T}}
\def\cU{{\cal U}} \def\bU{{\bf U}}  \def\mU{{\mathscr U}}
\def\cV{{\cal V}} \def\bV{{\bf V}}  \def\mV{{\mathscr V}}
\def\cW{{\cal W}} \def\bW{{\bf W}}  \def\mW{{\mathscr W}}
\def\cX{{\cal X}} \def\bX{{\bf X}}  \def\mX{{\mathscr X}}
\def\cY{{\cal Y}} \def\bY{{\bf Y}}  \def\mY{{\mathscr Y}}
\def\cZ{{\cal Z}} \def\bZ{{\bf Z}}  \def\mZ{{\mathscr Z}}

\def\Z{{\Bbb Z}}
\def\K{{\Bbb K}}
\def\R{{\Bbb R}}
\def\C{{\Bbb C}}
\def\T{{\Bbb T}}
\def\N{{\Bbb N}}
\def\S{{\Bbb S}}
\def\H{{\Bbb H}}
\def\dD{{\Bbb D}}

\def\qqq{\qquad}
\def\qq{\quad}
\newcommand{\ca}{\begin{cases}}
\newcommand{\ac}{\end{cases}}
\def\ma{\left(\begin{array}{cc}}    \def\am{\end{array}\right)}
\def\iint{\int\!\!\!\int}
\def\lt{\biggl}                     \def\rt{\biggr}
\let\ge\geqslant                   \let\le\leqslant
\def\[{\begin{equation}}            \def\]{\end{equation}}
\def\wt{\widetilde}                 \def\pa{\partial}
\def\sm{\setminus}                  \def\es{\emptyset}
\def\no{\noindent}                  \def\ol{\overline}
\def\iy{\infty}                     \def\ev{\equiv}
\def\/{\over}
\def\we{\wedge}
\def\ts{\times}
\def\os{\oplus}
\def\ss{\subset}
\def\h{\hat}
\def\wh{\widehat}
\def\Ra{\Rightarrow}
\def\ra{\rightarrow}
\def\la{\leftarrow}
\def\da{\downarrow}
\def\ua{\uparrow}
\def\lra{\leftrightarrow}
\def\Lra{\Leftrightarrow}
\def\Re{\mathop{\rm Re}\nolimits}
\def\Im{\mathop{\rm Im}\nolimits}
\def\supp{\mathop{\rm supp}\nolimits}
\def\sign{\mathop{\rm sign}\nolimits}
\def\Ran{\mathop{\rm Ran}\nolimits}
\def\Ker{\mathop{\rm Ker}\nolimits}
\def\Tr{\mathop{\rm Tr}\nolimits}
\def\const{\mathop{\rm const}\nolimits}
\def\Wr{\mathop{\rm Wr}\nolimits}
\def\diag{\mathop{\rm diag}\nolimits}
\def\dist{\mathop{\rm dist}\nolimits}

\def\wh{\widehat}  

\def\th{\theta}
\def\dlint{\displaystyle\int\limits}
\def\iintt{\mathop{\int\!\!\int\!\!\dots\!\!\int}\limits}
\def\intt{\mathop{\int\int}\limits}
\def\lim{\mathop{\rm lim}\limits}
\def\mult{\!\cdot\!}
\def\BBox{\hspace{1mm}\vrule height6pt width5.5pt depth0pt \hspace{6pt}}
\def\1{1\!\!1}
\newcommand{\bwt}[1]{{\mathop{#1}\limits^{{}_{\,\bf{\sim}}}}\vphantom{#1}}
\newcommand{\bhat}[1]{{\mathop{#1}\limits^{{}_{\,\bf{\wedge}}}}\vphantom{#1}}
\newcommand{\bcheck}[1]{{\mathop{#1}\limits^{{}_{\,\bf{\vee}}}}\vphantom{#1}}
\def\nh{\bhat}
\def\nc{\bcheck}
\newcommand{\oo}[1]{{\mathop{#1}\limits^{\,\circ}}\vphantom{#1}}
\newcommand{\po}[1]{{\mathop{#1}\limits^{\phantom{\circ}}}\vphantom{#1}}
\def\ctg{\mathop{\rm ctg}\nolimits}
\def\notto{\to\!\!\!\!\!\!\!/\,\,\,}

\def\pgbrk{\pagebreak}
\def\Dis{\mathop{\rm Dis}\nolimits}
\def\BBox{\hspace{1mm}\vrule height6pt width5.5pt depth0pt \hspace{6pt}}

%%%%%%%%%%%%%%%%%%%%%%%%%%%%%%%%%%%%%%%%%5
%% stary
%%%%%%%%%%%%%%%%%%%%%%%%%%%

\def\Twelve{
\font\Tenmsa=msam10 scaled 1200
\font\Sevenmsa=msam7 scaled 1200
\font\Fivemsa=msam5 scaled 1200
\textfont\msbfam=\Tenmsb
\scriptfont\msbfam=\Sevenmsb
\scriptscriptfont\msbfam=\Fivemsb

\font\Teneufm=eufm10 scaled 1200
\font\Seveneufm=eufm7 scaled 1200
\font\Fiveeufm=eufm5 scaled 1200
%\newfam\eufmfam
\textfont\eufmfam=\Teneufm
\scriptfont\eufmfam=\Seveneufm
\scriptscriptfont\eufmfam=\Fiveeufm}

\def\Ten{
\textfont\msafam=\tenmsa
\scriptfont\msafam=\sevenmsa
\scriptscriptfont\msafam=\fivemsa

\textfont\msbfam=\tenmsb
\scriptfont\msbfam=\sevenmsb
\scriptscriptfont\msbfam=\fivemsb

\textfont\eufmfam=\teneufm
\scriptfont\eufmfam=\seveneufm
\scriptscriptfont\eufmfam=\fiveeufm}

\title {Spectral estimates for matrix-valued periodic Dirac operators}

\author{
Evgeny Korotyaev
\begin{footnote}
{Institut f\"ur  Mathematik,  Humboldt Universit\"at zu Berlin,
Rudower Chaussee 25, 12489, Berlin, Germany,
e-mail: evgeny@math.hu-berlin.de
 }
\end{footnote}
}

\maketitle

\begin{abstract}

\no We consider the first order periodic systems perturbed by 
a $2N\ts 2N$ matrix-valued  periodic  potential on the real line.
The spectrum
of this operator is absolutely continuous and consists of
intervals separated by gaps. We define the Lyapunov function,
which is analytic on an associated N-sheeted Riemann surface. On each sheet the Lyapunov function has the standard properties of the Lyapunov
function for the scalar case. The Lyapunov function has branch points, which we call resonances. We prove the existence
of real or complex  resonances.
We determine the asymptotics of the periodic, anti-periodic spectrum and of the
resonances at high energy (in terms of the Fourier coefficients of the potential). We show that there exist two types of
gaps: i) stable gaps, i.e., the endpoints are periodic and
anti-periodic eigenvalues, ii) unstable (resonance) gaps, i.e., the
endpoints are resonances (real branch points). Moreover, we   determine various new trace formulae for potentials  and the  Lyapunov exponent.

\end{abstract}

%  \vskip 0.25cm


\section {Introduction and main results}
\setcounter{equation}{0}

Consider the self-adjoint operator $\cK$ acting on the space 
$L^2(\R)^{2N}$ and  given by
$$
\cK y=-iJ_1y'+V_ty, \ \ \ \ \ J_1=\ma
\!\!I_N\!\!&\!\! 0\!\!\\ \!\!0\!\!&\!\! -I_N\!\!\am,\ \ \ 
 V=\ma \!\!0\!\!&\!\! v\!\!\\
\!\! v^*&\!\!0\!\!\am,\ \ N\ge 1,
$$
here and below we use the notation
$(')=\pa /\pa t$ and $I_N$ is the identity $N\ts N$ matrix; 
$v$ is the {\bf complex} 1-periodic $N\ts N$  matrix and $V=V_t$ belongs to the Hilbert space $\mH$ given by
$$
\mH=\lt\{V=\ma \!\!0\!\!&\!\!v\!\!\\
\!\!v^*&\!\!0\!\!\am,
v=v^\top\!\!=\{v_{jk}(t)\}_{j,k=1}^N,\ \ t\in \R/\Z ,\ \
\|V\|^2=\!\!\int_0^1\! \Tr V_t^2dt<\iy 
\rt\}.
$$
Note that $\Re v$ and $\Im v$ are self-adjoint.
% equipped with the norm
%$\|V\|^2=\!\!\int_0^1\! \Tr V^2(t)dt<\iy $.
{\bf Without loss of generality  we assume} 
\[
\lb{12} \mV=\int_0^1V_t^2dt=\mV_0\os \mV_0,\ \ 
\mV_0={\rm diag} \{\n_1,...,\n_N\},\ \ \ \ 0\le \n_1\le \n_2\le...\le \n_N,
\]
here $\n_{1},..,\n_{N}$ are the eigenvalues of 
$\mV$, see Sect. 2 for the proof.

It is well known (see [DS] p.1486-1494, [Ge]) that the spectrum
$\s(\cK)$ of $\cK$ is absolutely continuous and consists of
non-degenerated intervals $\s_n,n\in \Z$.  These intervals are
separated by the gaps $g_n=(z_n^-,z_n^+)$ with the length
$|g_n|> 0,N_g^-< n< N_g^+$, where $-\iy\le N_g^-< N_g^+\le \iy$
and $N_g=N_g^+- N_g^--1$ is a total number of the gaps. 

Introduce the fundamental $2N\ts 2N$-matrix solutions
$\p(t,z)$ of the equation
\[
\lb{11}
-iJ_1\p'+V_t\p=z \p,\ \ \ \ z\in\C,\ \ \ \p(0,z)=I_{2N},
\]
and the monodromy $2N\ts 2N$-matrix
$\p(1,z)$. The matrix valued function $\p(1,\cdot)$ is entire. An eigenvalue of $\p(1,z)$ is called a {\it multiplier} of $\cK$: to each of them corresponds a solution $f$ of $-iJ_1f'+V_tf=z f$ with $f(t+1)=\t(z)f(t), t\in [0,1)$. They are
roots of the algebraic equation $D(\t,z)=0, \t,z\in\C$, where 
$
%\lb{ded}
D(\t,z)=\det (\p(1,z)-\t I_{N}),\  \t,z\in\C.
$  
Zeros of the function $\det (\cM(z)-I_N)$ (or $\det (\cM(z)-I_N)$) are periodic (or anti-periodic )  eigenvalues. 

There exist many papers about the first order periodic systems $N\ge 2$: Gel'fand and Lidskii [GL],Gohbert and Krein [GK], Krein [Kr], Potapov [Po], [YS]) and we mention new papers of Gesztesy and coauthors [CG],[CHGL], [GKM]. The basic results for direct spectral theory for the matrix case were obtained by Lyapunov and Poincar\'e  (see  [GL],[Kr],[YS]).

\no {\bf Theorem (Lyapunov, Poincar\'e)}. {\it For each $(V,z)\in \mH\ts \C$ 
the matrix-function $\p(1,z)$ satisfies:
\[
\lb{TL-1}
\p^{-1}(1,\cdot )=-J\p^\top (1,\cdot )J,
\]
\[
\lb{TL-2}
D(\t,\cdot ) =\t^{2N}D(\t^{-1},\cdot ),\ \ \ any \ \  \t\neq 0,\ \  where\ \ 
 D(\t,\cdot )=\det (\p(1,\cdot)-\t I_{2N}).
\]
\[
\lb{TL-3}
\s(\cK)=\{z\in \C: |\t(z)|=1 \qq{\rm   for\ some \ multiplier}
\ \t(z)\ \ {\rm of }\qq \p(1,z)\}.
\]
If for some $z\in\C$ (or $z\in\R$) $\t(z)$ is a multiplier
of multiplicity $d\ge 1$, then $\t^{-1}(z)$ (or $\ol\t(z)$) 
is a multiplier of multiplicity $d$. Moreover, each 
$\p(1,z), z\in\C$, has exactly $2N$ multipliers
$\t_j^{\pm 1}(z), j=1,..,N$. If $\t(z)$ is a simple multiplier and $|\t(z)|=1$, then
$\t'(z)\ne 0$.
}

The eigenvalues of $\p(1,z)$ are the zeros of the  equation $D(\t,z)=0$. This is an algebraic equation in $\t$ of degree $2N$ with coefficients, which are entire in $z\in\C$. It is well known
(see e.g. [Fo]) that the zeros $\t_j(z),j=1,..,2N$ 
of $D(\t,z)=0$ are (some branches of) analytic functions of $z$
with only algebraic singularities: the zeros $\t_j(z), j=1,..,2N$
constitute one or several branches of one or several analytic functions that have only algebraic singularities in $\C$.
Thus the number of eigenvalues of $\p(1,z)$ is a constant $N_e$
with the exception of some special values of $z$ 
(see below the definition of a resonance).
In general, there is a infinite number of such points
on the plane. If the functions $\t_j(z),j=1,..,2N$ are all distinct, then $N_e=2N$. If some of them are identical,
then we get $N_e<2N$ and $\p(1,z)$ is permanently degenerate.

Introduce the matrix-valued function
$\mL(z)={1\/2}(\p(1,z)+\p^{-1}(1,z)), z\in \C$
and the function $\F(z,\n)=\det (\mL(z)-\n I_{2N}), z,\n\in \C$.
Each zero of  $\F(\n,z)$ has multiplicity $\ge 2$ and define the Lyapunov function  by 
$
\D_j(z)={1\/2}(\t_j(z)+\t_j^{-1}(z)), \ j=1,..,N.
$
The Riemann surface for the multipliers $\t_j(z),j\in \N_{N}=\{1,..,N\}$  has $2N$ sheets, see \er{TL-2}. If $N=1$, then it has 2 sheets and the Lyapunov function is entire. Similarly, in the case $N\ge 2$ it is more convenient for us to  construct the Riemann surface for the Lyapunov function, which has N sheets. We need the following results from [K4].

\begin{theorem} \lb{T1}
Let $V\in \mH$. Then there exists an analytic function $\wt\D_s, s=1,..,N_0\le N$ on the
$N_s$-sheeted Riemann surface $\mR_s, N_s\ge 1$ having the following properties:

\no i) There exist disjoint subsets $\o_s,
s=1,..,N_0, \bigcup \o_s=\N_N$ such that all branches of
$\wt\D_s,s=1,2,..,N_0$ has the form
$\D_j(z)={1\/2}(\t_j(z)+\t_j^{-1}(z)), \ j\in \o_s$. Moreover,
 for any $z,\t\in \C$ the following relations hold true:
\[
\lb{T1-2} \det (\mL(z)-\n I_{2N})=\prod_1^{N_0} \F_s^2(\n,z),\
\ \ \ \ \F_s(\n,z)=\prod_{j\in \o_s}(\n-\D_j(z)),\ \ 
z,\n \in \C,
\]
\[
\lb{T1-2} \D_j(z)=\cos z+o(e^{|\Im z|}) \qq as \qqq |z|\to \iy,
\]
where the functions
$\F_s(\n,z)$ are entire with respect to $\n,z\in\C$. Moreover, if $\D_i=\D_j$ for
some $i\in \o_k, j\in \o_s$, then $\F_k=\F_s$ and $\wt \D_k=\wt
\D_s$.

\no ii) (The monotonicity property). Let some $\D_j,j=1,..,N
$, be real analytic on some interval
$Y=(\a,\b)\ss\R$ and $-1<\D_j(z)<1$ for any $z\in Y$.
Then $\D_j'(z)\ne 0$ for each $z\in Y$.

\no iii) The functions $\r,\r_s$ given by \er{T1-4} are entire,
\[
\lb{T1-4}
\r=\prod_{1}^{N_0}\r_s,\ \ \ 
\r_s(\cdot)=\!\!\!\!\prod_{i<j, i,j\in \o_s}\!\!\!\! (\D_i(\cdot)-\D_j(\cdot))^2.
\]
\no iv) The following identity holds true
\[
\s(\cK)=\R\sm \cup_{j=1}^N g_n,\qq g_n=(z_n^-,z_n^+), N_g^-<n<N_g^+
\]
where each gap $g_n=(z_n^-,z_n^+)$ is a bounded interval and $z_n^\pm$ are either periodic (anti-periodic) eigenvalues or real branch points of $\D_j$ (for some $j=1,..,N$)
which are zero of $\r$ (below we call such point a resonance).
\end{theorem}

%%%%%%%%%%%%%%%%%%%%%%%%%%%%%%%%%%%%%%%%%%%%%%%%%%%%%%%%%%%%%%%%%%

\no {\bf Remark.} 1) In the case of  $2\ts2$ system the
monodromy matrix has exactly 2 eigenvalues $\t,\t^{-1}$. The
Lyapunov function ${1\/2}(\t+\t^{-1})$  is an entire function of
the spectral parameter. It defines the band-gap structure of the
spectrum. By Theorem 1.1, the Lyapunov function for $2N\ts
2N$-matrix operator $\cK$ also defines the band-gap structure of
the spectrum, but it is the N-sheeted analytic function.

2) We have the following asymptotics (see Sect. 3)
\[
\lb{aD2} 
\D_j(z)=\cos z+{\sin z\/2z}\n_{j}+O\lt({e^{|\Im
z|}\/z^2}\rt), \ \ \  \ if \ \  V'\in \mH
\]
as $|z|\to \iy,\ j\in\N_N$.
Then firstly, $\r$ is not a polynomial since $\r$ is bounded on 
$\R$. Secondly, if $\n_{j'}\ne \n_{j}, j'\ne j$, then \er{aD2} implies  $\D_{j'}\ne \D_j$.

\no 3) In the case $N=2$ we determine $\D_1, \D_2, \r$ in terms of
the traces of the monodromy matrix.
Using \er{TL-2}, we have $D(\t,\cdot)=\t^4-T_1\t^3+{1\/2}(T_1^2-T_2)\t^2-T_1\t+1$,
which yields
\[
\lb{2l}
 D(\t,\cdot)=\lt(\t^2-2\D_1\t+1\rt)
\lt(\t^2-2\D_2\t+1\rt),\qq 
\D_1={T_1\/2}+{\sqrt{\r}\/2},\qqq \D_2={T_1\/2}-{\sqrt{\r}\/2},
\]
see [BBK], where $\r={T_2+4\/2}-{T_1^2\/4}$ and $T_m=\Tr \p(m,z),\ m=1,2$.

{\bf Definition.} {\it A zero $z_0\in\C$ of $\r$ given by \er{T1-4}
 is a {\bf resonance} of $\cK$}.

The main goal of this paper is to describe the spectrum of $\cK$
and to  determine the asymptotics of gaps and 
 resonances, periodic and anti-periodic eigenvalues  at high energy.
 We show that all resonances are real at high energy.
Moreover, we prove the existence of complex resonances
for some specific periodic potential. 
We have to underline that in the case of large $N$ the resonances
create the gaps in the spectrum of periodic operators, see  Theorem 
\ref{T2}  and remark after  Theorem \ref{T3}. If $N=1$, then 2-periodic eigenvalues create the gaps in the spectrum. 
In the present  paper we use some techniques from [K4] and [BBK], [BK], [CK].

 The periodic eigenvalues ($n$ is even) satisfy
\[
\lb{per}
..\le z_{N}^{-2,+}\le \underbrace{z_{1}^{0,+}\le z_{2}^{0+}\le ...\le z_{N}^{0,+}}_{n=0}\le 
\underbrace{z_{1}^{2,-}\le z_{1}^{2,+}\le
...\le z_{N}^{2,-}\le z_{N}^{2,+}}_{n=2}\le z_{1}^{4,-}\le \dots ,
\]
the anti-periodic eigenvalues ($n$ is odd) satisfy
\[
\lb{ape} 
..\le z_{N}^{-1,+}\le  \underbrace{z_{1}^{1,-}\le z_{1}^{1,+}\le
...\le z_{N}^{1,-}\le z_{N}^{1,+}}_{n=1}\le \underbrace{z_{1}^{3,-}\le z_{1}^{3,+}\le  ...
\le z_{N}^{3,+}}_{n=3}\le z_{1}^{5,-}\le \dots
\]
and they have asymptotics
\[
\lb{aspe} 
z_{j}^{n,\pm}=\pi n+o(1)\ \ \ as\ \ \  n\to \pm \iy,\ \ \ j\in \N_N=\{1,2,..,N\}.
\]
If $V=0$, then  these eigenvalues have the form $z_{j}^{n,\pm}=\pi n, (n,j)\in \Z\ts \N_N$.

Let $\wh {V_{n}'}=\!\!\int_0^1\!\! V_t'e^{i2\pi ntJ_1}dt,\ \  $
We formulate our first main result

\begin{theorem}   \lb{T2}
 Let $V,V'\in \mH$ and let $\z_{j}^{n,\pm}, j\in \N_N$ be
eigenvalues of the matrix $\mV-iJ_1\wh {V_n'}$. Then
 the periodic and anti-periodic eigenvalues
have the following asymptotics:
\[
\lb{T2-1} z_{j}^{n,\pm}=\pi n+{\z_{j}^{n,\pm}\/2\pi n}+O(n^{-2}), \ \ \
j\in \N_N\qq as \qq n\to \pm \iy.
\]
Assume that $\n_j\ne \n_{j'}$ for all $j\ne j'\in \o_s$ for some 
$s=1,..,N_0$. Then the function $\r_s$ has the zeros
$z_\a^{n\pm}, \a=(j,j'), j< j', j,j'\in \o_s, n\in \Z$, which are real at large $|n|$ and satisfy
\[
\lb{T2-2}
z_\a^{n\pm}=\pi n+{\n_{j}+\n_{j'}\/4 \pi n}
+O\lt({|\wh {V'_{n}}|\/n}+{1\/n^2}\lt),\ \ \ \alpha=(j,j')\ \ {\rm as}\ \  n\to \pm\iy.
\]

Let in addition $\n_{1}<..<\n_{N}$. Then for each $s=1,..,N_0$
and for large $n\to \pm\iy$ there exists a system of real intervals (gaps) $ g_\a^n=(z_\a^{n-},z_\a^{n+})$ such that 
$$
z_{j,j'}^{n\pm}=z_{j',j}^{n\pm},\ \ \ \ \a=(j,j'),\  j,j'\in \o_s,
\ \ \ \ \ z_{j,j_1}^{n-}\le z_{j,j_1}^{n+}<z_{j,j_2}^{n-}\le z_{j,j_2}^{n+}<...< z_{j,j_{N_s}}^{n-}\le z_{j,j_{N_s}}^{n+}, \ 
$$
\[
\lb{T2-3}
(-1)^n\D_j(z)>1, z\in g_{j,j}^n,\ \ and \ \  \ \ \ 
\ol \D_j(z)=\D_{j'}(z), z\in g_{j,j'}^n, \ \ {\rm if } \ \ j\ne j'
\]
i) Each branch $\D_j$ is real and is analytic on the set $(\pi n-{\pi\/2},\pi n+{\pi\/2})\sm\cup_{p\ne j} g_{p,j}^n$ and is not real on $\cup_{p\ne j} g_{p,j}^n$.

\no ii) If $z_\a^{n-}\ne z_\a^{n+}$ for some $\a=(j,j'), j\ne j'$, then  $z_\a^{n\pm}$ is the simple branch point (resonance) for the functions $\D_j, \D_{j'}$. If $z_\a^{n-}=z_\a^{n+}$, then $\D_j, \D_{j'}$ are analytic at $z_\a^{n\pm}$. 

\no iii) The following asymptotics hold true:
\[
\lb{T2-4} z_\a^{n\pm}=\pi n+{\n_{j}+\n_{j'}\pm |\wh {v_{n,\a}'}|\/2\pi n}+O\lt(|\wh {v_{n,\a}'}|+{1\/n}\lt),\ \ \
\wh {v_{n,\a}'}=\!\!\int_0^1\!\! v_\a'(t)e^{-i2\pi nt}dt,\ \  
\a=(j,j').
\]
\end{theorem}

{\bf Remark.} 1) $N_g=N_g^+-N_g^--1$ is the total number of gaps in the spectrum of $\cK$.

\no 2) If  $\n_{1}<..<\n_{N}$, then there exist infinite number of
resonances $z_\a^{n\pm}, \a=(j,j'), j\ne j'$ which form the gaps
in the spectrum of $\cK$, see below Theorem \ref{T3}.
Roughly speaking, resonances "form" the gaps, the number
of periodic and anti-periodic eigenvalues is less than 
the number of resonances. Thus there exists big difference
between $N=1$ and large $N$. In the first case
the endpoints of the gaps are 2-periodic eigenvalues.
In the second case, roughly speaking, the endpoints of the gaps are  resonances.

In the second main result we describe finite band potentials.

\begin{theorem}
\lb{T3} 
Let $V,V'\in \mH$ and let $\n_{1}<...<\n_{N}$.\\
(i) If the identity $\n_1+\n_N=\n_2+\n_{N-1}=...=\n_N+\n_1$ is not fulfilled, then $|N_g^\pm |<\iy$.\\
(ii) If $\n_1+\n_N=...=\n_N+\n_1$ holds true and there exists a sequence of indicies
$n_k\to\pm\iy$ such that $|\wh {v_{n,\a}'}|^2\!+|\!n_k|^{-1}=o(|\wh {v_{n,\a}'}|)$ as $k\to\pm\iy$, for each $\a=(j,N+1-j),  j\in \N_N$, then $N_g^\pm=\pm\iy$.

\end{theorem}

\no {\bf Remark.} 1) Consider $v=\diag\{v_{1},v_{2},..,v_{N}\}$,
i.e., the case when "variables are separated".
The transformation
$\mU_0: y=(y_1,..,y_{2N})^\top\to (y_1,y_{N+1},y_2,y_{N+2},...,y_N,y_{2N})^\top$ gives
$$
\mU_0\cK \mU_0^*=\os_1^N\cK_j,\ \ \ \ \ \cK_j=-i{\bf j_1}
 {d\/dt}+V_j,\qq {\bf j_1}=\ma \!\! 1\!\!&\!\! 0 \!\!\\ \!\! 0 &-1\am
\ \ \ \ \ V_j=\ma \!\!0\!\!&\!\!v_j\!\!\\
\!\!v_j^*&\!\!0\!\!\am.
$$ 
The operator $\cK_j$ for the case $N=1$
is well studied [YS], [K1-3], [Mi]. We have $\n_j=\int_0^1|v_j(t)|^2dt\ge 0$.
If $\n_j<\n_p$ for some $j<p$, then the number of gaps is $N_g<\iy$.\\
2) Note that the condition $|\wh {v_{n,\a}'}|^2\!+\!|n|^{-1}=o(|\wh {v_{n,\a}'}|)$,
$\a=(j,N+1-j),j\in \N_N$ as $n\to\pm\iy$, holds true for "generic"\ potentials $V, V'\in \mH$. This yields the existence
of the real resonance gaps $(z_\a^{n-},z_\a^{n+})$
at high energy. The coefficients $\wh {v_{n,\a}'},\a=(j,N+1-j), j\in \N_N$ 
(the second diagonal of the matrix $v$) "create" the gaps.

{\bf Example of complex resonances.} 
Consider the operator $\cK_{\n,\t}=-iJ_1{d\/dt}+V_{\n,\t},\n=1, {1\/2}, {1\/3},.., \t\in \R$ acting in $L^2(\R)^4$, where the real periodic potential $V_{\n,\t}$ is given by
\[
\lb{Vt} V_{\n,\t}=\ma \!\!0\!\!&\!\! v_{\n,\t}\!\!\\
\!\! v_{\n,\t}&\!\!0\!\!\am,\ \
v_{\n,\t}=-\ma a&\t b_\n\\\t b_\n&0\am,\qq {a\/2\pi}\in \R_+\sm \N,\qq
b_\n\in C(\T),
\]
%$\n=1, {1\/2}, {1\/3},..,.$ where  $\int_0^1|b_\n(t)|dt=1$
\[
\lb{Vn}
\int_0^1|b_\n(t)|dt=1,\qq
 \int_0^mb_\n(t)f(t)dt\to \int_0^m\d_{per}(t)f(t)dt,\qq
% {\rm as} \ \ \n\to 0,\ \ \
\d_{per}(t)\ev\sum_{-\iy}^\iy\d(t-n-{1\/2}),
\]
 as $\n\to 0$, for any $f\in C(0,m), m\in \N$.

 If $\t=0$, then the operator $\cK_{\n,0}=-iJ_1{d\/dt}+V_{\n,0}$ has the constant potential $V_{\n,0}$. In this case
there are no gaps in the spectrum of $\cK_{\n,0}$ and all resonances are given by $r_{n}^0=\pi n+{a^2\/4\pi n},n\in \Z\sm \{0\}$ with multiplicity
2. We show that there exist the non-degenerated resonance gaps for
small  $\t,\n$. In this example some resonances are real and
some are complex.

\begin{proposition}   \lb{T5}
{\it Let a potential $V_{\t,\n}$ satisfy \er{Vt}, \er{Vn}. 
 Then for each large integer  $n_0\ge 1+a$  there exist sufficiently small $\n,\ve>0$ such that the following statements hold true: 

\no i) Each function $\r(z,V_{\t,\n}), \t\in (-\ve,\ve)$ in the disk
$z\in \dD_{n_0}=\{|z|<\pi n_0 +1\}$ has exactly $4n_0$ simple
zeros $r_{n,\n}^\pm(\t), 1\le \pm n\le n_0$,
where $r_{n,\n}^\pm(\t)$ is analytic function in the disk $\t\in \{|\t|<\ve\}$ and $r_{n,\n}^\pm(0)=r_n^0$ and the
following estimates hold
\[
\lb{p1} 
r_{n,\n}^\pm(\t)=r_n^0\pm \t (\sqrt{R_n}+o(1))\qq as \qq \t\to 0,\qq 
 where \qq 
\cases
{R_n<0  &if\ \ \ \ $|n|<{ a\/2\pi}$\cr
           R_n>0 &if\ \ \ \ $|n|>{ a\/2\pi}$\cr},
%\qq as \qq \t\to 0,
\]
 for some constant $R_n$ and if $\t \in (-\ve,\ve), |n|>{ a\/2\pi}$, then $(r_{n,\n}^-(\t),r_{n,\n}^+(\t))\ss \R$ is  a gap.
%There are no other roots in the disk $\dD_{n_0}$.

\no ii) Each function $D(1,z,V_{\t,\n})D(-1,z,V_{\t,\n})$ in the disk
$\dD_{n_0}$ has exactly $4n_0+4$ zeros
$z_{n,m}^\pm (\t,\n), \\ -n_0\le n\le n_0, m=1,2$ and the 
asymptotics $z_{n,m}^\pm (\t,\n)=z_{n,m}^{0,\pm}+o(1)$  holds as
$\t\to 0$.
 }
\end{proposition}

\no {\bf Remark.} 1) If $0<a<2\pi $, then $\r(\cdot,V_{\t,\n})$
has only real roots $r_{n,\n}^\pm(\t)$ in each large disk $\dD_{n_0}$
for sufficiently small $\t, \n$. If $1<{a\/2\pi
}$, then $\r(z,V_{\t,\n})$ has at least two non-real roots
$r_{1,\n}^\pm(\t)$ for small $\t, \n$.
2) We show that operator $\cK_{\t,\n}$ has 
new gaps (so-called resonance gaps). The endpoints of
the resonance gap are the branch points of the Lyapunov function,
and, in general, they are not the periodic (or anti-periodic)
eigenvalues. These endpoints are not stable. If they are real (see
\er{p1}), then we have a gap. If they are complex ($0< n\le {a\/2\pi}$), then we have not a gap, we have only the branch
points of the Lyapunov function in the complex plane.
3) We have a similar complicated distribution of other resonances, which
are poles of S-matrix for scattering  for Schr\"odinger operator with compactly supported potentials on the real line see [K5], [Z].

{\bf We consider the conformal mapping associated with the operator}
$\cK$. 
Introduce the simple conformal mapping $\e:\C\sm [-1,1]\to \{\z\in
\C: |\z|>1\}$ by
\[
\lb{dx} \e(z)=z+\sqrt{z^2-1},\ \ \ \ z\in \C\sm [-1,1],\ \ {\rm
and}\ \ \ \  \e(z)=2z+o(1)\qq as \qq |z|\to \iy.
\]
Note that $\e(z)=\ol \e(\ol z), z\in \C\sm [-1,1]$ since $\e(z)>1$
for any $z>1$. Due to the properties of the Lyapunov functions we have $|\e(\wt\D_s(\z))|>1, \z\in \mR_s^+=\{\z\in\mR_s: \Im \z>0\}$. Thus we can introduce the quasimomentum $k_j, j=1,2,..,N$ (we fix some branch of $\arccos$ and $\D_j(z)$) and the function $q_j$ by
\[
\lb{dkm} k_j(z)=\arccos \D_j(z)=i\log \e(\D_j(z)),\ \  
q_j(z)=\Im k_j(z)=\log |\e(\D_j(z))|,\ \ 
\]
$z\in \mR_0^+=\C_+\sm \b_+,
\b_+=\!\bigcup_{\b\in \cB_\D\cap \C_+}
[\b, \b+i\iy)$
where $\cB_\D$ is the set of all branch points of the function $\D$.
The branch points of $k_m$ belong to $\cB_\D$. 
Define the {\bf averaged quasimomentum} $k$, the {\bf
density} $p$ and the {\bf Lyapunov exponent} $q$ by
\[
k(z)=p(z)+iq(z)={1\/N}\sum_1^N k_j(z), \ \ \  q(z)=\Im k(z),
\ \  z\in \mR_0^+.
\]
For the function $k(z)=p(z)+iq(z),z\in \ol\C_+$ we introduce
formally integrals
\[
\lb{dQSI} Q_n={1\/\pi}\int_{\R}t^nq(t)dt,\ \ \ \ \
 I_n^S={1\/\pi}\int_{\R}t^nq(t)dp(t),\ \ \ \ 
 I_n^D={1\/\pi}\iint_{\C_+}|k_{(n)}'(z)|^2dxdy,
\]
$n=0,1,2,$ here and below $ k_{(0)}=k-z,  k_{(1)}=zk_{(0)}, z=x+iy\in\C$. 
 Let $C_{us}$ denote the class of all
real upper semi-continuous functions  $h:\R\to \R$. With any $h\in
C_{us}$  we associate the "upper" domain $\K(h)=\{k=p+iq\in\C:
q>h(p), p\in \R\}$. We formulate our last result.

\begin{theorem}   \lb{T4}
i) Let $V\in \mH$. Then the averaged quasimomentum $k={1\/N}\sum_1^{N} k_j$ is analytic
in $\C_+$ and $k:\C_+\to k(\C_+)=\K(h)$ is a conformal mapping for
some $h\in C_{us}$. Furthermore, for some branches $k_j, j=1..,N$ the following asymptotics, identities and estimates hold true:
 \[
 \lb{kc-1}
k(z)-z=-{Q_0+o(1)\/z}\qq as \qq |z|\to \iy, \  if\ y>r|x|, \ \ for \ any \ r>0,
\]
 \[
 \lb{kc-2}
Q_0=I_0^D+I_0^S={\|V\|^2\/4N},
\]
$$ 
q|_{\s_{(N)}}=0,\ \ \ 0< q^2|_{\s_{(1)}\cup g}\le 2Q_0,\ \ \  where
\ \ \ 
 \s_{(N)}\!=\!\{z\in\R: \D_1(z),..,\D_N(z)\in [-1,1]\},
$$
$$
\s_{(1)}\!=\!\{z\in \R: \D_m(z)\in [-1,1],\ \D_p(z)\notin [-1,1] \ {\rm \ some} \ m,p=1,..,N\},
$$ 
\[
 \lb{T3-4}
\sum_n |g_n|^2\le {4\|V\|^2\/N}. 
\]

\no ii) Let aditionly  $V'\in \mH$. Then  the following asymptotics,
identities  hold true:
 \[
 \lb{T4-1}
k(z)-z=-{Q_0\/z}-{Q_1\/z^2}-{Q_2+o(1)\/z^3}\qq as \qq |z|\to \iy, \  \ \  if\ y>r|x|, \ \ for \ any \ r>0,
\]
\[
 \lb{T4-2}  Q_1={\Tr\/8N}\int_0^1-iJ_1V_t'V_tdt,\  \ \ \ \ \ 
Q_2=I_1^D+I_2^S-{Q_0^2\/2}={\Tr\/16N}\int_0^1\rt({V_t'}^2+V_t^4\rt)dt.
\]
\end{theorem}

A priori estimates for various parameters of the Dirac operator
 (the norm of a periodic potential, effective masses,
gap lengths, height of slits, action variables for NLS and so on) were obtained in [KK1-2],  [K1-3], [Mi] only for the case $N=1$.  In order to get the required estimates  the authors of [KK1-2],[K1-3], [Mi],... used the global quasi-momentum as the conformal mapping, which was introduced into the spectral  theory of  the  Hill operator by Marchenko-Ostrovski [MO].

\begin{figure}[htb]
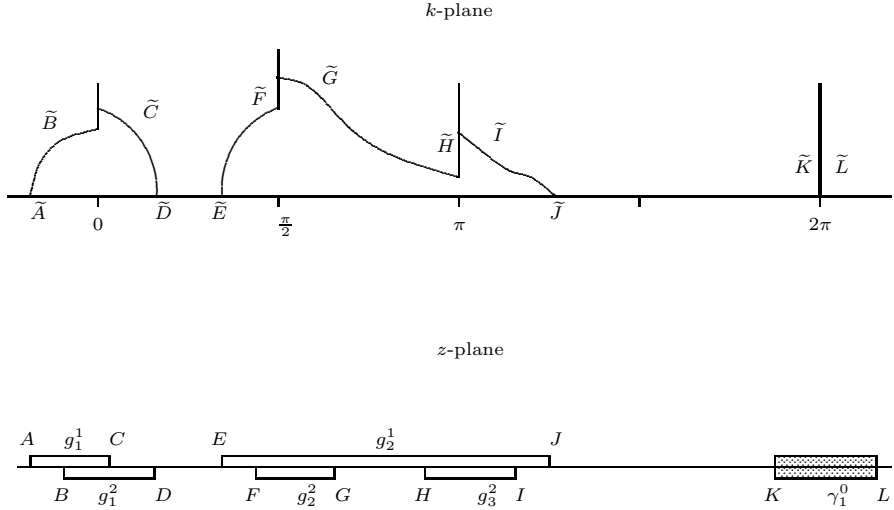

\begin{center}
\font\thinlinefont=cmr5
\begingroup\makeatletter\ifx\SetFigFont\undefined%
\gdef\SetFigFont#1#2#3#4#5{%
  \reset@font\fontsize{#1}{#2pt}%
  \fontfamily{#3}\fontseries{#4}\fontshape{#5}%
  \selectfont}%
\fi\endgroup%
\mbox{\beginpicture \setcoordinatesystem units
<0.62992cm,0.62992cm> \unitlength=0.62992cm \linethickness=1pt
\setplotsymbol ({\makebox(0,0)[l]{\tencirc\symbol{'160}}})
\setshadesymbol ({\thinlinefont .}) \setlinear
\linethickness=
0.500pt \setplotsymbol ({\thinlinefont .}) {\circulararc 72.861
degrees from  6.638 22.824 center at  7.300 21.129
}%
%
% Fig POLYLINE object
%
\linethickness= 0.500pt \setplotsymbol ({\thinlinefont .})
{\setshadegrid span <1pt> \shaderectangleson \putrectangle corners
at 17.145 15.479 and 19.287 15.240 \setshadegrid span <5pt>
\shaderectanglesoff
}%
%
% Fig POLYLINE object
%
\linethickness= 0.500pt \setplotsymbol ({\thinlinefont .})
{\setshadegrid span <1pt> \shaderectangleson \putrectangle corners
at 17.145 15.240 and 19.287 15.003 \setshadegrid span <5pt>
\shaderectanglesoff
}%
%
% Fig POLYLINE object
%
\linethickness= 0.500pt \setplotsymbol ({\thinlinefont .})
{\putrectangle corners at  5.476 15.479 and 12.383 15.240
}%
%
% Fig POLYLINE object
%
\linethickness= 0.500pt \setplotsymbol ({\thinlinefont .})
{\putrectangle corners at  9.762 15.240 and 11.667 15.003
}%
%
% Fig POLYLINE object
%
\linethickness= 0.500pt \setplotsymbol ({\thinlinefont .})
{\putrectangle corners at  6.191 15.240 and  7.857 15.003
}%
%
% Fig POLYLINE object
%
\linethickness= 0.500pt \setplotsymbol ({\thinlinefont .})
{\putrectangle corners at  2.142 15.240 and  4.047 15.003
}%
%
% Fig POLYLINE object
%
\linethickness= 0.500pt \setplotsymbol ({\thinlinefont .})
{\putrectangle corners at  1.429 15.479 and  3.095 15.240
}%
%
% Fig POLYLINE object
%
\linethickness= 0.500pt \setplotsymbol ({\thinlinefont .})
{\putrule from  1.156 15.240 to 19.6 15.240
}%
%
% Fig POLYLINE object
%
\linethickness= 0.500pt \setplotsymbol ({\thinlinefont .})
{\putrule from  0.953 20.955 to 19.6 20.955
}%
%
% Fig POLYLINE object
%
\linethickness= 0.500pt \setplotsymbol ({\thinlinefont .})
{\putrule from  6.638 23.467 to  6.644 23.467 \plot  6.644 23.467
6.659 23.465 / \plot  6.659 23.465  6.687 23.461 / \plot  6.687
23.461  6.723 23.457 / \plot  6.723 23.457  6.769 23.451 / \plot
6.769 23.451  6.822 23.442 / \plot  6.822 23.442  6.879 23.434 /
\plot  6.879 23.434  6.936 23.423 / \plot  6.936 23.423  6.993
23.410 / \plot  6.993 23.410  7.046 23.398 / \plot  7.046 23.398
7.095 23.383 / \plot  7.095 23.383  7.140 23.368 / \plot  7.140
23.368  7.184 23.351 / \plot  7.184 23.351  7.224 23.332 / \plot
7.224 23.332  7.262 23.309 / \plot  7.262 23.309  7.300 23.285 /
\plot  7.300 23.285  7.338 23.258 / \plot  7.338 23.258  7.368
23.235 / \plot  7.368 23.235  7.400 23.207 / \plot  7.400 23.207
7.434 23.180 / \plot  7.434 23.180  7.468 23.148 / \plot  7.468
23.148  7.501 23.116 / \plot  7.501 23.116  7.537 23.080 / \plot
7.537 23.080  7.573 23.044 / \plot  7.573 23.044  7.612 23.004 /
\plot  7.612 23.004  7.650 22.964 / \plot  7.650 22.964  7.688
22.921 / \plot  7.688 22.921  7.728 22.879 / \plot  7.728 22.879
7.766 22.837 / \plot  7.766 22.837  7.806 22.792 / \plot  7.806
22.792  7.846 22.748 / \plot  7.846 22.748  7.885 22.705 / \plot
7.885 22.705  7.925 22.663 / \plot  7.925 22.663  7.963 22.623 /
\plot  7.963 22.623  7.999 22.583 / \plot  7.999 22.583  8.037
22.543 / \plot  8.037 22.543  8.073 22.507 / \plot  8.073 22.507
8.109 22.471 / \plot  8.109 22.471  8.145 22.435 / \plot  8.145
22.435  8.181 22.401 / \plot  8.181 22.401  8.217 22.367 / \plot
8.217 22.367  8.255 22.333 / \plot  8.255 22.333  8.293 22.301 /
\plot  8.293 22.301  8.333 22.267 / \plot  8.333 22.267  8.376
22.236 / \plot  8.376 22.236  8.418 22.204 / \plot  8.418 22.204
8.460 22.172 / \plot  8.460 22.172  8.507 22.140 / \plot  8.507
22.140  8.551 22.109 / \plot  8.551 22.109  8.600 22.077 / \plot
8.600 22.077  8.647 22.047 / \plot  8.647 22.047  8.695 22.020 /
\plot  8.695 22.020  8.744 21.990 / \plot  8.744 21.990  8.791
21.965 / \plot  8.791 21.965  8.839 21.937 / \plot  8.839 21.937
8.888 21.912 / \plot  8.888 21.912  8.937 21.888 / \plot  8.937
21.888  8.983 21.865 / \plot  8.983 21.865  9.032 21.844 / \plot
9.032 21.844  9.078 21.823 / \plot  9.078 21.823  9.127 21.802 /
\plot  9.127 21.802  9.172 21.785 / \plot  9.172 21.785  9.218
21.766 / \plot  9.218 21.766  9.267 21.747 / \plot  9.267 21.747
9.318 21.728 / \plot  9.318 21.728  9.370 21.709 / \plot  9.370
21.709  9.426 21.689 / \plot  9.426 21.689  9.487 21.668 / \plot
9.487 21.668  9.553 21.647 / \plot  9.553 21.647  9.622 21.624 /
\plot  9.622 21.624  9.696 21.601 / \plot  9.696 21.601  9.777
21.575 / \plot  9.777 21.575  9.859 21.550 / \plot  9.859 21.550
9.946 21.524 / \plot  9.946 21.524 10.031 21.497 / \plot 10.031
21.497 10.116 21.471 / \plot 10.116 21.471 10.196 21.448 / \plot
10.196 21.448 10.270 21.425 / \plot 10.270 21.425 10.331 21.408 /
\plot 10.331 21.408 10.382 21.391 / \plot 10.382 21.391 10.420
21.380 / \plot 10.420 21.380 10.446 21.374 / \plot 10.446 21.374
10.458 21.370 / \plot 10.458 21.370 10.465 21.368 /
}%
%
% Fig POLYLINE object
%
\linethickness= 0.500pt \setplotsymbol ({\thinlinefont .})
{\putrule from  6.668 22.80 to  6.668 24.052
}%
%
% Fig POLYLINE object
%
\linethickness=1pt \setplotsymbol
({\makebox(0,0)[l]{\tencirc\symbol{'160}}}) {\putrule from 18.098
20.955 to 18.098 23.336
}%
%
% Fig POLYLINE object
%
\linethickness= 0.500pt \setplotsymbol ({\thinlinefont .}) {\plot
10.465 22.316 10.471 22.310 / \plot 10.471 22.310 10.486 22.299 /
\plot 10.486 22.299 10.511 22.278 / \plot 10.511 22.278 10.547
22.250 / \plot 10.547 22.250 10.592 22.214 / \plot 10.592 22.214
10.645 22.172 / \plot 10.645 22.172 10.698 22.130 / \plot 10.698
22.130 10.753 22.085 / \plot 10.753 22.085 10.806 22.043 / \plot
10.806 22.043 10.854 22.005 / \plot 10.854 22.005 10.899 21.969 /
\plot 10.899 21.969 10.941 21.937 / \plot 10.941 21.937 10.977
21.907 / \plot 10.977 21.907 11.013 21.880 / \plot 11.013 21.880
11.045 21.855 / \plot 11.045 21.855 11.074 21.831 / \plot 11.074
21.831 11.102 21.808 / \plot 11.102 21.808 11.134 21.783 / \plot
11.134 21.783 11.165 21.759 / \plot 11.165 21.759 11.197 21.734 /
\plot 11.197 21.734 11.229 21.711 / \plot 11.229 21.711 11.261
21.685 / \plot 11.261 21.685 11.292 21.662 / \plot 11.292 21.662
11.324 21.639 / \plot 11.324 21.639 11.356 21.618 / \plot 11.356
21.618 11.388 21.596 / \plot 11.388 21.596 11.417 21.577 / \plot
11.417 21.577 11.447 21.558 / \plot 11.447 21.558 11.474 21.541 /
\plot 11.474 21.541 11.502 21.526 / \plot 11.502 21.526 11.527
21.514 / \plot 11.527 21.514 11.553 21.503 / \plot 11.553 21.503
11.576 21.493 / \plot 11.576 21.493 11.604 21.482 / \plot 11.604
21.482 11.633 21.471 / \plot 11.633 21.471 11.663 21.463 / \plot
11.663 21.463 11.692 21.457 / \plot 11.692 21.457 11.722 21.448 /
\plot 11.722 21.448 11.756 21.440 / \plot 11.756 21.440 11.788
21.433 / \plot 11.788 21.433 11.822 21.425 / \plot 11.822 21.425
11.853 21.414 / \plot 11.853 21.414 11.887 21.406 / \plot 11.887
21.406 11.919 21.395 / \plot 11.919 21.395 11.949 21.385 / \plot
11.949 21.385 11.980 21.372 / \plot 11.980 21.372 12.012 21.357 /
\plot 12.012 21.357 12.040 21.342 / \plot 12.040 21.342 12.067
21.325 / \plot 12.067 21.325 12.097 21.306 / \plot 12.097 21.306
12.131 21.283 / \plot 12.131 21.283 12.167 21.258 / \plot 12.167
21.258 12.207 21.226 / \plot 12.207 21.226 12.251 21.192 / \plot
12.251 21.192 12.298 21.154 / \plot 12.298 21.154 12.349 21.114 /
\plot 12.349 21.114 12.397 21.074 / \plot 12.397 21.074 12.442
21.035 / \plot 12.442 21.035 12.480 21.004 / \plot 12.480 21.004
12.507 20.983 / \plot 12.507 20.983 12.522 20.968 / \plot 12.522
20.968 12.531 20.961 /
}%
%
% Fig POLYLINE object
%
\linethickness= 0.500pt \setplotsymbol ({\thinlinefont .})
{\putrule from 10.478 21.368 to 10.478 23.336
}%
%
% Fig POLYLINE object
%
\linethickness= 0.500pt \setplotsymbol ({\thinlinefont .})
{\putrule from  2.857 22.384 to  2.855 22.384 \plot  2.855 22.384
2.849 22.382 / \plot  2.849 22.382  2.832 22.377 / \plot  2.832
22.377  2.800 22.371 / \plot  2.800 22.371  2.758 22.360 / \plot
2.758 22.360  2.703 22.348 / \plot  2.703 22.348  2.639 22.333 /
\plot  2.639 22.333  2.574 22.316 / \plot  2.574 22.316  2.508
22.299 / \plot  2.508 22.299  2.445 22.280 / \plot  2.445 22.280
2.385 22.265 / \plot  2.385 22.265  2.330 22.248 / \plot  2.330
22.248  2.282 22.231 / \plot  2.282 22.231  2.237 22.214 / \plot
2.237 22.214  2.197 22.200 / \plot  2.197 22.200  2.159 22.183 /
\plot  2.159 22.183  2.125 22.166 / \plot  2.125 22.166  2.091
22.147 / \plot  2.091 22.147  2.057 22.128 / \plot  2.057 22.128
2.026 22.106 / \plot  2.026 22.106  1.994 22.083 / \plot  1.994
22.083  1.962 22.060 / \plot  1.962 22.060  1.930 22.035 / \plot
1.930 22.035  1.901 22.007 / \plot  1.901 22.007  1.869 21.977 /
\plot  1.869 21.977  1.839 21.946 / \plot  1.839 21.946  1.812
21.914 / \plot  1.812 21.914  1.782 21.882 / \plot  1.782 21.882
1.757 21.848 / \plot  1.757 21.848  1.731 21.816 / \plot  1.731
21.816  1.708 21.783 / \plot  1.708 21.783  1.687 21.749 / \plot
1.687 21.749  1.666 21.717 / \plot  1.666 21.717  1.649 21.685 /
\plot  1.649 21.685  1.632 21.654 / \plot  1.632 21.654  1.615
21.620 / \plot  1.615 21.620  1.600 21.588 / \plot  1.600 21.588
1.587 21.556 / \plot  1.587 21.556  1.573 21.520 / \plot  1.573
21.520  1.560 21.484 / \plot  1.560 21.484  1.547 21.444 / \plot
1.547 21.444  1.535 21.402 / \plot  1.535 21.402  1.522 21.355 /
\plot  1.522 21.355  1.509 21.304 / \plot  1.509 21.304  1.494
21.249 / \plot  1.494 21.249  1.482 21.194 / \plot  1.482 21.194
1.469 21.137 / \plot  1.469 21.137  1.456 21.086 / \plot  1.456
21.086  1.446 21.040 / \plot  1.446 21.040  1.439 21.002 / \plot
1.439 21.002  1.433 20.976 / \plot  1.433 20.976  1.431 20.961 /
\plot  1.431 20.961  1.429 20.955 /
}%
%
% Fig POLYLINE object
%
\linethickness= 0.500pt \setplotsymbol ({\thinlinefont .})
{\putrule from  2.857 22.384 to  2.857 23.336
}%
%
% Fig POLYLINE object
%
\linethickness= 0.500pt \setplotsymbol ({\thinlinefont .})
{\putrule from 18.098 20.955 to 18.098 20.743
}%
%
% Fig POLYLINE object
%
\linethickness= 0.500pt \setplotsymbol ({\thinlinefont .})
{\putrule from 14.287 20.955 to 14.287 20.743
}%
%
% Fig POLYLINE object
%
\linethickness= 0.500pt \setplotsymbol ({\thinlinefont .})
{\putrule from 10.478 20.955 to 10.478 20.743
}%
%
% Fig POLYLINE object
%
\linethickness= 0.500pt \setplotsymbol ({\thinlinefont .})
{\putrule from  6.668 20.955 to  6.668 20.743
}%
%
% Fig POLYLINE object
%
\linethickness= 0.500pt \setplotsymbol ({\thinlinefont .})
{\putrule from  2.857 20.955 to  2.857 20.743
}%
%
% Fig TEXT object
%
\put{\SetFigFont{7}{8.4}{\rmdefault}{\mddefault}{\updefault}{$\wt L$}%
} [lB] at 18.404 21.465
%
% Fig TEXT object
%
\put{\SetFigFont{7}{8.4}{\rmdefault}{\mddefault}{\updefault}{$\wt K$}%
} [lB] at 17.554 21.465
%
% Fig TEXT object
%
\put{\SetFigFont{7}{8.4}{\rmdefault}{\mddefault}{\updefault}{$g_1^1$}%
} [lB] at  2.142 15.708
%
% Fig TEXT object
%
\put{\SetFigFont{7}{8.4}{\rmdefault}{\mddefault}{\updefault}{$k$-plane}%
} [lB] at  9.762 24.765
%
% Fig TEXT object
%
\put{\SetFigFont{7}{8.4}{\rmdefault}{\mddefault}{\updefault}{$z$-plane}%
} [lB] at 10.001 17.621
%
% Fig TEXT object
%
\put{\SetFigFont{7}{8.4}{\rmdefault}{\mddefault}{\updefault}
{$\gamma_1^0$}%
} [lB] at 18.098 14.527
%
% Fig TEXT object
%
\put{\SetFigFont{7}{8.4}{\rmdefault}{\mddefault}{\updefault}
{$g_2^1$}%
} [lB] at  8.572 15.716
%
% Fig TEXT object
%
\put{\SetFigFont{7}{8.4}{\rmdefault}{\mddefault}{\updefault}
{$g^2_3$}%
} [lB] at 10.715 14.527
%
% Fig TEXT object
%
\put{\SetFigFont{7}{8.4}{\rmdefault}{\mddefault}{\updefault}
{$g_2^2$}%
} [lB] at  6.905 14.527
%
% Fig TEXT object
%
\put{\SetFigFont{7}{8.4}{\rmdefault}{\mddefault}{\updefault}{$g_1^2$}%
} [lB] at  2.857 14.527
%
% Fig TEXT object
%
\put{\SetFigFont{7}{8.4}{\rmdefault}{\mddefault}{\updefault}{$J$}%
} [lB] at 12.383 15.716
%
% Fig TEXT object
%
\put{\SetFigFont{7}{8.4}{\rmdefault}{\mddefault}{\updefault}{$I$}%
} [lB] at 11.667 14.527
%
% Fig TEXT object
%
\put{\SetFigFont{7}{8.4}{\rmdefault}{\mddefault}{\updefault}{$H$}%
} [lB] at  9.525 14.527
%
% Fig TEXT object
%
\put{\SetFigFont{7}{8.4}{\rmdefault}{\mddefault}{\updefault}{$G$}%
} [lB] at  7.857 14.527
%
% Fig TEXT object
%
\put{\SetFigFont{7}{8.4}{\rmdefault}{\mddefault}{\updefault}{$F$}%
} [lB] at  5.952 14.527
%
% Fig TEXT object
%
\put{\SetFigFont{7}{8.4}{\rmdefault}{\mddefault}{\updefault}{$D$}%
} [lB] at  4.047 14.527
%
% Fig TEXT object
%
\put{\SetFigFont{7}{8.4}{\rmdefault}{\mddefault}{\updefault}{$B$}%
} [lB] at  1.905 14.527
%
% Fig TEXT object
%
\put{\SetFigFont{7}{8.4}{\rmdefault}{\mddefault}{\updefault}{$A$}%
} [lB] at  1.190 15.716
%
% Fig TEXT object
%
\put{\SetFigFont{7}{8.4}{\rmdefault}{\mddefault}{\updefault}{$C$}%
} [lB] at  3.095 15.716
%
% Fig TEXT object
%
\put{\SetFigFont{7}{8.4}{\rmdefault}{\mddefault}{\updefault}{$E$}%
} [lB] at  5.239 15.716
%
% Fig TEXT object
%
\put{\SetFigFont{7}{8.4}{\rmdefault}{\mddefault}{\updefault}{$L$}%
} [lB] at 19.287 14.527
%
% Fig TEXT object
%
\put{\SetFigFont{7}{8.4}{\rmdefault}{\mddefault}{\updefault}{$K$}%
} [lB] at 16.906 14.527
%
% Fig TEXT object
%
\put{\SetFigFont{7}{8.4}{\rmdefault}{\mddefault}{\updefault}{$\wt J$}%
} [lB] at 12.383 20.479
%
% Fig TEXT object
%
\put{\SetFigFont{7}{8.4}{\rmdefault}{\mddefault}{\updefault}{$\wt I$}%
} [lB] at 11.191 22.147
%
% Fig TEXT object
%
\put{\SetFigFont{7}{8.4}{\rmdefault}{\mddefault}{\updefault}{$\wt H$}%
} [lB] at 10.001 21.907
%
% Fig TEXT object
%
\put{\SetFigFont{7}{8.4}{\rmdefault}{\mddefault}{\updefault}{$\wt G$}%
} [lB] at  7.586 23.338
%
% Fig TEXT object
%
\put{\SetFigFont{7}{8.4}{\rmdefault}{\mddefault}{\updefault}{$\wt F$}%
} [lB] at  6.088 22.9
\put{\SetFigFont{7}{8.4}{\rmdefault}{\mddefault}{\updefault}{$\wt E$}%
} [lB] at  5.239 20.479
%
% Fig TEXT object
%
\put{\SetFigFont{7}{8.4}{\rmdefault}{\mddefault}{\updefault}{$\wt C$}%
} [lB] at  3.810 22.623
%
% Fig TEXT object
%
\put{\SetFigFont{7}{8.4}{\rmdefault}{\mddefault}{\updefault}{$\wt B$}%
} [lB] at  1.666 22.384
%
% Fig TEXT object
%
\put{\SetFigFont{7}{8.4}{\rmdefault}{\mddefault}{\updefault}{$\wt D$}%
} [lB] at  4.047 20.479
%
% Fig TEXT object
%
\put{\SetFigFont{7}{8.4}{\rmdefault}{\mddefault}{\updefault}{$\wt A$}%
} [lB] at  1.429 20.479
%
% Fig TEXT object
%
\put{\SetFigFont{7}{8.4}{\rmdefault}{\mddefault}{\updefault}{$0$}%
} [lB] at  2.75 20.253
%
% Fig TEXT object
%
\put{\SetFigFont{7}{8.4}{\rmdefault}{\mddefault}{\updefault}
{$\frac{\pi}{2}$}%
} [lB] at  6.478 20.253
%
% Fig TEXT object
%
\put{\SetFigFont{7}{8.4}{\rmdefault}{\mddefault}{\updefault}{$\pi$}%
} [lB] at 10.35 20.253
%
% Fig TEXT object
%
\put{\SetFigFont{7}{8.4}{\rmdefault}{\mddefault}{\updefault}{$2\pi$}%
} [lB] at 17.858 20.253
% Fig CIRCULAR ARC object
%
\linethickness= 0.500pt \setplotsymbol ({\thinlinefont .})
{\circulararc 76.531 degrees from  4.098 20.961 center at  2.290
21.099
}%
\linethickness=0pt \putrectangle corners at  0.906 25.146 and
21.002 14.385
\endpicture}

\end{center}
\caption{The domain $\K(h)=k(\C_+)$ and the graph of the function $h$
} \label{fig1}
\end{figure}

The mapping
$k:\C_+\to \K(h)$ is illustrated in Figure \ref{fig1} for the case
$N=2$. The
integral $I_0^S\ge 0$ is the area between the boundary of $\K(h)$
and the real line. In Figure 1 the upper picture is a domain
$\K(h)$ and the points $\wt A=k(A), \wt B=k(B),...$  . The
spectral interval $(A,B)$ (with multiplicity 2) of the z-domain
is mapped on the curve $(\wt A,\wt B)$ of the k-domain, the
interval (a gap)  $(B,C)$ of the z-domain  is mapped on the
vertical slits, which lies on the line $\Re k=0$. The spectral
interval $(C,D)$ (with multiplicity 2) of the z-domain is mapped
on the curve $(\wt C,\wt D)$ of the k-domain. The spectral
interval $(D,E)$ (with multiplicity 4) of the z-domain is mapped
on the interval $(\wt D,\wt C)$ of the k-domain. The case of the
interval $(E,J)$ is similar. The resonace gap $(K,L)$  of the
z-domain is mapped on the vertical slits, which lies on the line
$\Re k=2\pi$. In fact we have the graph of the function $h(p),
p\in \R$, which coinsides with the boundary of $\K(h)$.

%%%%%%%%%%%%%%%%%%%%%%%%%%%%%%       FIGURE 1  %%%%%%%%%%%%%%
%%%%%%%%%%%%%%%%%

%%%%%%%%%%%%%%%%%%%%%%%%%%%%%%   END  OF  FIGURE 1  %%%%%%%%%%%%%%%%%

We describe the plan of our paper. In Sect. 2 we obtain the basic properties of fundamental solution $\p(t,z)$. In Sect. 3 we determine the asymptotics of the fundamental solution $\p(t,z)$ for the case $V, V'\in \mH$.  In Sect. 4 we determine the asymptotics of the Lyapunov function and multipliers at high energy and prove Theorem \ref{T2}, \ref{T4}. In Sect. 5 we prove Proposition \ref{T5}. In Sect. 5 we 
determine the asymptotics of $\det \mL(z)$ as $\Im z\to \iy$.

\section {The fundamental solutions}
\setcounter{equation}{0}

In this section we study $\p$. We begin with some notational
convention. A vector $h=\{h_n\}_1^N\in \C^N$ has the Euclidean
norm $|h|^2=\sum_1^N|h_n|^2$, while a $N\ts N$ matrix $A$ has the
operator norm given by $|A|=\sup_{|h|=1} |Ah|$. Note that
$|A|^2\le \Tr A^*A$. 

Recall the identity: if $A$ is $N\ts N$ matrix, then  
\[
\lb{poD}
\det (A-\n I_N)=(-1)^N\sum_0^{N}a_j\n^{N-j},\qq
  a_0=1,\qq a_1=-\cA_1,\qq \f_2=-{\cA_2+\cA_1\f_1\/2},
  \]
and $.., \qq \f_j=-{1\/j}\sum_1^{j}\cA_ka_{j-k},..,\qq a_N=\det A, \qq \cA_m(z)=\Tr A^m$.

Below we need the identity
\[
\lb{21} J_1V=-VJ_1,  \ \ \  e^{zJ_1} V=Ve^{-zJ_1},\  
\]
for any $z\in \C$. The solution of the equation \er{11} satisfies the integral equation
\[
\lb{22} \p (t,z)=\p_0(t,z)-i\int_0^te^{izJ_1(t-s)}J_1V_s\p (s,z)ds,\ \
\ \p_0(t,z)=e^{iztJ_1},\ \ t\ge 0,\ z\in \C.
\]
It is clear that Eq. \er{22}  has a solution as a power series in
$V$ given by
\[
\lb{23} \p (t,z)=\sum_{n\ge0}\p_n(t,z), \ \ \
\p_n(t,z)=-i\int_0^te^{izJ_1(t-s)}J_1V_s\p_{n-1}(s,z)ds,\ \ n\ge1.
\]
 Using \er{21}, \er{23} we have
\[
\lb{24} \p_1(t,z)=-i\int_0^te^{izJ_1(t-s)}J_1V_se^{izsJ_1}ds=
-i\int_0^te^{izJ_1(t-2s)}J_1V_sds,
\]
\[
\lb{25} \p_2=-i\int_0^te^{izJ_1(t-t_1)}J_1V_{t_1}\p_{1}(t_1,z)dt_1=
\int_0^tdt_1\int _0^{t_1}e^{izJ_1(t-2t_1+2t_2}V_{t_1}V_{t_2}dt_2.
\]
Proceeding by induction, we obtain
%\[
%\lb{26} \p_{n}(t,z)=(-i)^n \int_0^tdt_1\int _0^{t_1}dt_2\dots \int
%_0^{t_{n-1}} e^{izJ_1(t-2t_1+2t_2\dots +(-1)^n2t_n)}J_1V_{t_1}\dots
%J_1V_{t_n}dt_n,\]and
\[
\lb{26} \p_{2n}(t,z)=\int_0^tdt_1\dots \int _0^{t_{2n-1}}
e^{izJ_1(t-2t_1+2t_2\dots +2t_{2n})}V_{t_1}\dots V_{t_{2n}}dt_{2n},
\]
\[
\lb{27} \p_{2n+1}(t,z)=-i\int_0^xdt_1\dots \int _0^{t_{2n}}
e^{izJ_1(t-t_1+2t_2\dots -2t_{2n+1})}J_1V_{t_1}\dots
V_{t_{2n+1}}dt_{2n+1}.
\]

We need the following results from [K4].

\no  \begin{lemma}\lb{T21} 
Let $V\in\mH$.  For each $z\in \C$ there
exists a unique solution $\p$ of Eq. \er{22} given by \er{23} and
series \er{23} converge uniformly on bounded subsets of
$\R\ts\C\ts \mH$. For each $t\ge 0$ the function $\p (t,z)$ is
entire on $\C$.  Moreover, for any $n\ge 0$ and $(t,z)\in
[0,\iy)\ts \C$ the following estimates and asymptotics hold true:
\[
\lb{213} |\p_n(t,z)|\le {e^{|\Im z|t}\/n!}\lt(\int_0^t |V_s|ds
\rt)^n,
\]
\[
\lb{214} |\p(t,z)-\sum_0^{n-1}\p_j(t,z)|\le {(\sqrt
t\|V\|)^n\/n!}e^{t|\Im z|+\int_0^t |V_s|ds},
\]
\[
\lb{28} \p(t,z)-e^{iztJ_1}=o(e^{t|\Im z|})\ \ {\rm as} \ \ 
|z|\to \iy,
\] 
 uniformly on bounded $t\in \R$.
 If the sequence $V^{\n}\to V$ weakly in $\mH$, as $\n\to \iy$, then $\p (t,z,V^{\n})\to \p (t,z,V)$
 uniformly on bounded subsets of $\R\ts \C$.
\end{lemma}

Below we need the simple properties of matrices
$a, b, c\in \cA=\lt\{A=\ma 0&a_1\\a_2&0\am :\ a_1,a_2 \ {\rm is \ }\ N\ts N \ {\rm matrix} \lt\}$ given by
\[
\lb{29} abc, \  J_1a,\ \ e^{zJ_1}a\in \cA,\ \
{\rm all}\ \ \ z\in \C,
\]
\[
\lb{210} ab=(-a_1b_1+a_2b_2)+(a_1b_2-a_2b_1)J_1,
\]
\[
\lb{211} \Tr a=0,\ \   \Tr Ja^n=0,\ \  n\ge 0.
\]
For any matrices $A, B$ the following identities
hold
\[
\lb{212}
 \Tr AB=\Tr BA, \ \ \ \ \ \ \ \ \ol{\Tr A}=\Tr A^*.
\]

 Using \er{25}-\er{27} we define the function
 \[
\lb{217}
T_{j,1}(z,V)=
\Tr\int_0^jV_{s_1}\int_0^{s_1}V_{s_2}e^{izJ_1(j-2s_1+2s_2)}ds,
\]
\[
\lb{216}
T_{j,n}(z,V)=\Tr\p_{2n}(j,z)=\Tr\int_0^j\!\!\!\dots
\!\!\! \int_0^{t_{2n-1}}\!\!\!e^{izJ_1(j-2s_1+2s_2\dots
+2s_{2n})}V_{s_1}\dots V_{s_{2n}}ds,\ n\ge 2,
\]
where $s=(s_1,..,s_{2n})\in \R^{2n}$. 

\begin{lemma}  \lb{T22} 
Let $V\in\mH$. The functions $T_j(\cdot,V), j=1,2,..,N$
are entire on $\C$ and $T(z,V)\in \R$  for all $z\in \R$.
Moreover, the function $T(z,V)$ satisfies
\[
\lb{218} T_j(z,V)=T_j(z,-V)=2N\cos jz+\sum _{n\ge 1}T_{j,n}(z,V),
\]
\[
\lb{219} |T_j(z,V)|\le 2Ne^{ j(|\Im z|+\|V\|)},
\]
\[
\lb{220}
T_j(z,V)=2N\cos jz +o(e^{j|\Im z|})\qq as  \qqq |z|\to \iy,
\]
and  $\Tr \p _{2n+1}(t,z)=0$.
Series \er{218} converge uniformly on bounded subsets of $\C\ts
\mH$. If a sequence $V^{\n}$ converges weakly to $V$ in $\mH$
as $\n\to \iy$, then $T_j(z,V^{\n})\to T_j(z,V)$ uniformly on
bounded subsets of $\C$.
\end{lemma}
 \no {\bf Proof.}
By Lemma \ref{21}, series \er{218} converge uniformly and
absolutely on bounded subsets of $\C\ts \cH$. Each term in
\er{218} is an entire function of $z$, then $T_j$ is an entire
function of $z,V$.
Moreover, if the sequence $V^{\n}$ converges weakly to $V$ in
$\cH$, as $\n\to \iy$, then $T_j(z,V^{\n})\to T_j(z,V)$
uniformly on bounded subsets of $\C$.

 We have $T_j=\Tr \p(j,z)=\Tr \sum_{n\ge
0}\p_n(j,z)$ and $\Tr \p_0(j,z)=2N\cos jz$ for $j\ge 1$. The
estimate \er{213}  yields \er{219} and  \er{28} gives
\er{220}.

The first relation in \er{29} yields   $ V_{t_1}\dots
V_{t_{2n+1}}\in \cA$ for any $t_1,\dots t_{2n+1}\in\R$. Then
relations in \er{29} give
$\Tr e^{izJ_1(t-2t_1+2t_2\dots -2t_{2n+1})}J_1V_{t_1}\dots V_{t_{2n+1}}=0,
$
for any $t, t_1,\dots t_{2n+1}\in \R$, which together with
\er{28} implies $\Tr \p_{2n+1}(t,z)=0$.

Below  we will show that $T(z,V)\in \R$  for all $z\in \R$
(see \er{229}).
 $\BBox$

{\bf We will show} \er{12}.
Consider the self-adjoint operator $Jy'+\O (t)y, J=\ma
\!\!0\!\!&\!\! I_N\!\!\\ \!\!-I_N\!\!&\!\! 0\!\!\am$, acting on
the Hilbert space $\os_1^{2N} L^2(\R)$, where the real matrix $\O $ is given by
\[
\lb{221}
\O =\ma \!\!\O _1\!\!&\!\!\O _2\!\!\\
\!\!\O _2&\!\!-\O _1\!\!\am :
 \ \  \O_2=\O _2^*,\ \  \O _1=\O _1^*,\ \ \ 
\|\O \|^2=\!\!\int_0^1\!\!\!\!\Tr \O ^2(t)dt<\iy.
\]
Let $M(t,z), t\in \R$ be the fundamental solution of the equation
$JM'+\O M=zM, M(0,z)=I_{2N}$. Note that $M(t,z)$ is real for $z\in \R$.

Define the unitary matrix $\mU={1\/\sqrt 2}(J_1+iJ)=\mU^*, \mU^2=I_{2N}$.
Using the identities 
\[
\lb{222}
J_2=\ma
\!\!0\!\!&\!\! I_N\!\!\\ \!\!I_N\!\!&\!\! 0\!\!\am,\ \ 
\mU J\mU=-iJ_1,\ \ \ \mU J_1\mU=iJ,\ \ \ \mU J_2\mU=-J_2,\ \ \ 
\O=J_1\O_1+J_2\O_2,\ 
\]
we deduce that $M_c=\mU M\mU$ satisfy the equation
$-iJ_1M_c'+\O_cM_c=zM_c, M_c(0,z)=I_{2N}$, where $\O_c$ is given by
\[
\lb{223}
 \O _c=\O _c^*=\mU\O \mU=\ma \!\!0\!\!&\!\! \o\!\!\\
\!\! \o^*&\!\!0\!\!\am \in \mH,\ \ \o=-\O _2+i\O _1=
\o^\top .
\]
Thus we obtain
\[
\lb{224}
\int_0^1 \O _c^2(t)dt=\mV_{1}\os \mV_{2},
%\ma \mV_{1}& 0\\ 0& \mV_{2}\am
\ \ \mV_{1}=\int_0^1 \o(t)\o^*(t)dt,
\ \ \ \ \mV_{2}=\int_0^1 \o^*(t)\o(t)dt,
\]
and 
\[
\lb{225}
\mV_{1}=E \mV_{0} E^*=\ol \mV_{2}\ge 0,\ \ 
\ \ \mV_{0}=\diag \{\n_1,..,\n_N\},
\]
for some unitary matrix $E$ and the diagonal matrix $\mV_0$.
Define the unitary matrix $\mE=E\os \ol E$. The function $\p(t,z)=\mE^* M_c(t,z)\mE$ satisfies  
$-iJ_1\p'+V\p=z \p,\ \  \p(0,z)=I_{2N}$, where 
\[
\lb{227}
\mE=E\os \ol E,\ \ V=\mE^* \O_c\mE=\ma \!\!0\!\!&\!\! v\!\!\\
\!\! v^*&\!\!0\!\!\am,\ \ \ \ v=E^*\o\ol E,
\]
\[
\lb{228}
\int_0^1 V_t^2dt=\mV_0\os \mV_0,
\]
\[
\lb{229}
\Tr \p(t,z)=\Tr M(t,z),
\]
which gives $V\in \mH$ and $\Tr \p(1,z)\in \R$ for all $z\in \R$.

It is well known that for real $\O$ we have $M(t,z)JM(t,z)^\top=J$
(see [GL], [YS]). Then $M_c=\mU M\mU$ and \er{222} give
$$
-iJ_1=M_c(1,z)(-iJ_1)\mU \mU^\top (M_c(1,z)^\top)\mU^\top \mU,\ \ 
(-iJ_1)\mU\mU^\top =-J,
$$
which yields
\[
M_c(1,z)JM_c(1,z)^\top=J.
\]
The similar arguments and $\p=\mE^*M_c\mE$ imply
$$
\p\mE^*J\mE^{*\top}\p^\top \mE^\top=\mE^*J,\ \ \ \ \ 
\mE^*J\mE^{*\top}=J,
$$
which yields $\p J\p^\top=J$ and \er{TL-1} is proved for $V\in \mH$.

\section {Estimates of $\p$ for the case $V'\in \mH$}
\setcounter{equation}{0}

Sect. 2 does not give the needed estimates of the fundamental
solution $\p$ at high energy. In order to determine the
asymptotics of $\p$ we will do some modification. Define the
integral operator $K$ and the matrix-valued function $a_t(z)$ by
\[
\lb{31} (Kf)(t)=\int_0^t\!\!\!e_{t-s}W_sf(s)ds,\ \ \ 
W=-iJ_1V^2-V',\ e_t=e_t(z)=e^{itzJ_1},\ a_t(z)=I-{V_t\/2z} .
\]
where $e_t=e_t(z)=e^{itzJ_1}$. Introduce
\[
\lb{32} \wh {V'}(z)=\int_0^1\!\!\!V_t'e^{i2tzJ_1}dt,\ \ \ 
\wh {V^3}(z)=\int_0^t\!\!\!V_t^3e^{i2tzJ_1}dt,\ \ \ \mV=\int_0^1 V_t^2dt.
\]

\begin{lemma}  \lb{T31} 
For each $(z,V')\in (\C\sm \{0\})\ts\cH $
the solution $\p=a^{-1}\P a_0$, where $\P$ satisfies 
\[
\lb{33} \P=\p_0+\ve Ka^{-1}\P,\ \ \ \ \
\P=a\p a_0^{-1},\ \ a_0=I-\ve V_0,\ \ \ve={1\/2z},
\]
\[
\lb{34} \P=\p_0+\sum_{n\ge 1}\P_n,\ \ \ \P_n=\ve^n (Ka^{-1})^n\p_0,
\]
where series \er{34} converge uniformly on bounded subsets of
$\R\ts(\C\sm \{0\})\ts \cH$. Moreover, if $\sup_{t\in \R} |V_t|\le |z|$, then for any $j-1,m\in \N$ the following estimates are fulfilled:
\[
\lb{35} 
|\P_n(t,z)|\le {e^{|\Im z|t}\/n!|z|^n}\lt(\int_0^t
|W_s|ds\rt)^n,
\]
\[
\lb{36}
 |\P(m,z)-\p_0(m,z)-\sum_1^{j-1}\P_n(m,z)|\le
{\vk^j\/j!}e^{ m(|\Im z|+\vk)},\ \ \ \ \ \
\vk\ev{\|V\|^2+\|V'\|\/|z|},\ \
\]
\[
\lb{37}
 \P_1(1,z)=-{e^{izJ_1}\/2z}
 \rt(iJ_1\mV+\wh {V'}(z)+{1\/2z}\rt(iJ_1\wh {V^3}(z)+\int_0^1
V_t'V_tdt\rt)\rt)+O({e^{|\Im z|}\/z^3}),\ \ 
\]
\[
\lb{38}
 \P_2(1,z)=- {e^{izJ_1}\/4z^2}\!\!\int_0^1\!\!dt\!\!\int_0^t\!\!\rt(V_t^2V_s^2-iJ_1V_t^2V_s'e_{2s}-iV'_te_{2t}J_1V_s^2-V'_te_{2t}V_s'e_{2s}\rt)ds+O({e^{|\Im z|}\/z^3})
\]
as $|z|\to \iy$ and where $e_t=e^{iztJ_1}$.
\end{lemma}

 \no {\bf Proof.} 
Using \er{21}, \er{22}, $\ve={1\/2z}$ and integrating by parts
we get
$$
\p(t,z)-e^{iztJ_1}=-i\int_0^t J_1e^{izJ_1(t-s)}V_s\p(s,z)ds=
-i\int_0^t e^{izJ_1(t-2s)}J_1\rt(V_se^{-izsJ_1}\p(s,z)\rt)ds
$$
$$
=\ve e^{izJ_1(t-2s)}\rt(V_se^{-izsJ_1}\p(s,z)\rt)\rt|_0^t-
\ve\int_0^t e^{izJ(t-2s)}\rt(V_se^{-izsJ_1}\p(s,z)\rt)'ds
$$
$$
=\ve \rt(V_t\p(t,z)-e^{iztJ_1}V_0\rt)-\ve\int_0^t e^{izJ_1(t-2s)}\rt(V'_se^{-izsJ_1}+ie^{izsJ_1}J_1V_s^2\rt)\p(s,z)ds.
$$
Thus we obtain
$$
a(t,z)\p(t,z)=e^{iztJ}a(0,z)+\ve\int_0^t e^{izJ(t-s)}W_s\p(s,z)ds
$$
which yields \er{33}. We will show \er{34}-\er{36}. Using 
$|1/a(t,z)|\le 2$ for $\sup_{t\in \R} |V_t|\le |z|$ and $|\p_0(t,z)|\le e^{|\Im z|t}$, we have
$$
|\P_n(t,z)|\le 2|\ve| \int_0^t e^{|\Im
z|(t-t_1)}|W_{t_1}||\P_{n-1}(t_1,z)|dt_1
$$
$$
\le
(2|\ve|)^n\!\!\int_0^t\!\!dt_1\!\!\int_0^{t_1}\!\!\!\!dt_2\dots\!\!\!
\int _0^{t_{n-1}}\!\! e^{|\Im z|t}|W_{t_1}|\dots |W_{t_n}| dt_1
\le {(2|\ve|)^n\/n!}e^{|\Im z|t}\lt(\int_0^t |W_{t_1}|dt_1 \rt)^n,
$$
which gives \er{35}.
Estimates \er{35} and $\int_0^m|W(t)|dt\le m(\|V'\|+\|V\|^2)$ imply \er{34} and \er{36}. We get
\[
\P_1(1,z)=\ve Ka^{-1}\p_0=\ve K\p_0+\ve^2 KV\p_0+O(\ve^3e^{
|\Im z|}).
\]
Recall $e_t=e^{iztJ_1}$. \er{214} implies
$$
K\p_0=-\int_0^1e_{1-t}(iJ_1V_t^2+V_t')e_tdt=-ie_1J_1\mV-e_1\wh {V'}(z),
$$$$
KV\p_0=-\int_0^1e_{1-t}(iJ_1V_t^2+V_t')V_te_tdt=-ie_1J_1\wh {V^3}(z)-e_1\int_0^1V_t'V_tdt,
$$
which yields \er{37}. Consider the second term 
$\P_2=\ve^2 (Ka^{-1})^2\p_0=\ve^2 \P_{20}+O(\ve^3e^{|\Im z|})$, where
$$
\P_{20}=
\int_0^1\!\!dt\!\!\int_0^te_{1-t}(iJ_1V_t^2+V'_t)e_{t-s}
(iJ_1V_s^2+V_s')e_sds
$$
$$
=\int_0^1dt\int_0^te_{1-t}\rt[-J_1V_t^2e_{t-s}J_1V_s^2+iJ_1V_t^2e_{t-s}V_s'-iV'_te_{t-s}J_1V_s^2+ V'_te_{t-s}V_s'\rt]e_s
$$
$$
=-e_{1}\int_0^1dt\int_0^t\rt[V_t^2V_s^2-ie_{-2s}J_1V_t^2V_s'-ie_{-2t}V'_tJ_1V_s^2-e_{-2t+2s}V'_tV_s'\rt]ds
$$
which gives \er{38}.
$\BBox$

In order to determine the asymptotics of the Lyapunov function
we need the following modification. Substituting
$\p^{-1}=-J\p^\top J$ into $\P(1,z)^{-1}=a_0\p(1,z)^{-1}a_0^{-1}$ we get 
$$
\P(1,z)^{-1}=-a_0Ja_0\P(1,z)^\top a_0^{-1}Ja_0^{-1}=-cJ\P(1,z)^\top Jc^{-1},\ \ c=I-\ve^2V_0^2,
\ve={1\/2z},
$$
which  yields
\[
\lb{310}
L={1\/2}\rt(\P(1,z)+\P(1,z)^{-1}\rt)
={1\/2}\rt(\P(1,z)-cJ\P^\top(1,z)Jc^{-1}\rt).
\]
We determine the asymptotics of $\P(1,z)$. Using
\er{36}, we have
$$
\P=\p_0+\ve Ka^{-1}\p_0+\ve^2 Ka^{-1}Ka^{-1}\p_0+\ve^2 Ka^{-1}Ka^{-1}Ka^{-1}\p_0+O(\ve^4e^{|\Im z|})
$$
$$
=\p_0+\ve K(I+\ve V+\ve^2 V^2)\p_0+\ve^2 K(I+\ve V)K(I+\ve V)\p_0+\ve^2 K^3\p_0+O(\ve^4e^{|\Im z|})
$$
$$
=\p_0+\ve \vk_1+\ve^2 \vk_2+\ve^3\vk_3+O(\ve^4e^{|\Im z|}),
$$
where
\[
\lb{311}
\vk_1=K\p_0,\ \vk_2=(KV+K^2)\p_0,\ \ \vk_3=(KV^2+K^2V+KVK+K^3)\p_0.
\]
We get
$$
2L=\p_0+\ve \vk_1+\ve^2 \vk_2+\ve^3\vk_3-(I-\ve^2V_0^2)J(
\p_0+\ve \vk_1+\ve^2 \vk_2+\ve^3\vk_3)^\top J(I+\ve^2V_0^2)+
O(\ve^4e^{|\Im z|})
$$
$$
=(\p_0-J\p_0^\top J)+\ve (\vk_1-J\vk_1^\top J)+\ve^2(\vk_2-J\vk_2^\top J+V^2J\p_0^\top J-J\p_0^\top JV_0^2)
$$$$+
\ve^3(\vk_3-J\vk_3^\top J+V_0^2J\vk_1^\top J-J\vk_1^\top JV_0^2)+O(\ve^4e^{|\Im z|}).
$$
Thus we have
\[
\lb{312}
L=\cos z+{\ve}L_1+{\ve^2}L_2+\ve^3L_3
+O(\ve^4e^{|\Im z|}),
\]
where
\[
\lb{313}
L_1={\vk_1-J\vk_1^\top J\/2},\ \ \ \ L_2={\vk_2-J\vk_2^\top J\/2},
L_3={1\/2}(\vk_3-J\vk_3^\top J+V_0^2J\vk_1^\top J-J\vk_1^\top JV_0^2).
\]
Below we need the identities
\[
\lb{314}
J^*V^\top J=-V, \ \ \ \ J^*W^\top J=-W,\ \ \ \ \ J^*J_1J=-J_1.
\]
We determine the asymptotics of $L$. 
Let $u_t=\int_0^tV_s^2ds,\ f_t=\int_0^tV_s'e_{2s}ds$.

\begin{lemma}  \lb{T32} 
If $V, V'\in \mH$, then asymptotics \er{aD2}, the following identities and  asymptotics  are fulfilled:
\[
\lb{315} L_1=\mV\sin z, \ \ \ z\in \C,\ \ where\ \ 
\mV= \int_0^1V_t^2dt,
\]
\[
\lb{316} 
L_2=i{e_1J_1\/2}(\mV\wh {V'}+\wh {V'}\mV)-i\sin z\int_0^1J_1V_t'V_tdt+L_{21}+L_{22},
\]
$$
L_{21}=-{1\/2}\int_0^1\rt(e_{1}u'u+e_{-1}uu'\rt)dt,\ \ \ \
L_{22}={1\/2}\int_0^1\rt(e_{1}f'f+e_{-1}ff'\rt)dt,
$$
\[
\lb{317}  L_{22}=o(e^{|\Im z|}),
\]
\[
\lb{318} L_2(\pi n)={(-1)^{n}\/2}\rt(-\mV^2+iJ_1(\mV\wh {V_n'}+\wh {V_n'}\mV)+(\wh {V_n'})^2\rt),\ \ \wh {V_n'}=\wh {V'}(\pi n),
\] 
\[
\lb{319} 
\D_j(z)=\cos (z-{\n_j\/2z})+{O(|\wh {V_n'}|+|n|^{-1})\/n^2}, \ \ \ as \ z=\pi n+O\rt({1\/n}\rt), \ j=1,..,N.
\]
%where $[A,B]=AB-BA$ for matrix $A,B$.

\end{lemma}
\no {\bf Proof.} Recall $e_t=e^{iztJ_1}$. \er{311},\er{31}  give 
$$
J^*\vk_1^\top J=J^*(K\p_0)^\top J=J^*\int_0^1(e_{1-t}W_te_t)^\top Jdt=
\int_0^1e_{-t}J^*W_t^\top Je_{t-1}dt=-\int_0^1e_{-t}W_te_{t-1}dt.
$$
Then
$$
2L_1=\!\!\int_0^1\!\!\rt(e_{1-t}W_te_t-e_{-t}W_te_{t-1}\rt)dt=
-iJ_1\!\!\int_0^1\!\!\rt(e_{1-t}V_t^2e_t-e_{-t}V_t^2e_{t-1}\rt)dt=
-iJ_1\!\!\int_0^1\!\!(e_{1}-e_{-1})V_t^2dt
$$
which yields \er{315}. We determine $L_2={1\/2}(\vk_2-J\vk_2^\top J)$. Using \er{311},\er{31} we get
$$
\vk_2=(KV+K^2)\p_0=\int_0^1e_{1-t}W_t\rt(V_te_t+\int_0^t(e_{t-s}W_se_sds\rt)dt
$$
$$
J^*\vk_2^\top J=J^*\!\!\int_0^1\rt(e_{t}V_t^\top +\!\!\int_0^t(e_{s}W_s^\top e_{s-t}ds\rt)W_t^\top e_{1-t}dtJ
%$$$$
=\!\!\int_0^1\!\!\rt(e_{-t}V_t+\int_0^t(e_{-s}W_se_{s-t}ds\rt)W_te_{t-1}dt.
$$
This yields
$$
L_2=F+S,
%A_{21}={1\/2}\int_0^1\rt(e_{1-t}W_tV_te_t+e_{-t}V_tW_te_{t-1}\rt)%dt,
\ \ 
S={1\/2}\int_0^1dt\int_0^t\rt(e_{1-t}W_te_{t-s}W_se_s+
e_{-s}W_se_{s-t}W_te_{t-1}\rt)ds
$$
and using $W=-iJ_1V^2-V'$ we obtain
$$
F=\!\!\int_0^1\!\!\rt(e_{1-t}W_tV_te_t+e_{-t}V_tW_te_{t-1}\rt){dt\/2}
%$$$$
=-\!\!\int_0^1\!\!\rt(e_{1-t}(iJ_1V_t^2+V_t')V_te_t+e_{-t}V_t(iJ_1V_t^2+V_t')e_{t-1}\rt){dt\/2}
$$$$
=-{1\/2}\int_0^1\rt(e_{1}V_t'V_t+e_{-1}V_tV_t'\rt)dt
%$$$$
={e_{-1}-e_{1}\/2} \int_0^1V_t'V_tdt=
-iJ_1\sin z\int_0^1V_t'V_tdt
$$
since $\int_0^1V_tV_t'dt=-\int_0^1V_t'V_tdt$.  Consider the second term,
$$
2S=\int_0^1dt\int_0^t\rt(e_{1-t}(iJ_1V_t^2+V'_t)e_{t-s}
(iJ_1V_s^2+V_s')e_s+
e_{-s}(iJ_1V_s^2+V_s')e_{s-t}(iJ_1V_t^2+V'_t)e_{t-1}\rt)ds
$$
$$
=\int_0^1dt\int_0^t\rt(e_{1-t}\rt[-J_1V_t^2e_{t-s}J_1V_s^2+iJ_1V_t^2e_{t-s}V_s'+iV'_te_{t-s}J_1V_s^2+ V'_te_{t-s}V_s'\rt]e_s+
$$$$
+e_{-s}\rt[-J_1V_s^2e_{s-t}J_1V_t^2+iJ_1V_s^2e_{s-t}V_t'+iV_s'e_{s-t}J_1V_t^2
+V'_se_{s-t}V_t'\rt]e_{t-1}
$$$$
=\int_0^1dt\int_0^t\rt(\rt[-e_{1}V_t^2V_s^2+ie_{1-2s}J_1V_t^2V_s'+ie_{1-2t}V'_tJ_1V_s^2+e_{1-2t+2s}V'_tV_s'\rt]
$$
$$
+\rt[-e_{-1}V_s^2V_t^2+ie_{1-2t}J_1V_s^2V_t'+ie_{1-2s}V_s'J_1V_t^2
+e_{-1+2t-2s}V'_sV_t'\rt]\rt)ds
$$$$
=-\int_0^1dt\rt(e_{1}u'u+e_{-1}uu'\rt)dt+{ie_1J_1\/2}(\mV \wh {V'}+\wh {V'}\mV )+
\int_0^1dt\rt(e_{1}f'f+e_{-1}ff'\rt)dt
$$
where  $\wh {V'}(z)=\int_0^1V'(t)e_{2t}dt$ and $u=\int_0^tV_s^2ds,\ f=\int_0^tV_s'e_{2s}ds$, and here we used
$$
\int_0^1dt\int_0^t\rt(e_1V_t^2V_s^2+e_{-1}V_s^2V_t^2\rt)dtds=
\int_0^1dt\rt(e_{1}u'u+e_{-1}uu'\rt)dt,
$$
$$
\int_0^1(u'f+f'u+uf'+fu')dt=u_1f_1+f_1u_1=\mV \wh {V'}+\wh {V'}\mV,
$$
$$
\int_0^1dt\int_0^t\rt(e_{1-2t+2s}V'_tV_s'+e_{2t-2s-1}V'_sV_t'\rt)ds=
\int_0^1dt\rt(e_{1}f'f+e_{-1}ff'\rt)dt.
$$
We have $e_{\pm 1}=(-1)^nI_{2N}$ at $z=\pi n$ and then
\[
\lb{320} 
L_{21}(\pi n)=-{(-1)^n\/2}\int_0^1\rt(u'u+uu'\rt)dt
=-(-1)^n{\mV^2\/2},
\]
\[
\lb{321} 
L_{22}(\pi n)={(-1)^n\/2}\int_0^1\rt(f'f+ff'\rt)dt={(-1)^n\/2}f^2(\pi n)={(-1)^n\/2}(\wh {V'})^2(\pi n),
\]
which yields \er{318}. Using \er{312},\er{315}\er{318}
we obtain
\[
\lb{322}
L(z)=\cos (z-\ve \mV )+\ve^2\rt(L_2(\pi n)+(-1)^n{\mV^2\/2}\rt)+O(\ve^3)
\ \  as \ z=\pi n+O(1/n). 
\]
Recall the simple fact: Let $A, B$ be matrices and 
 and $\s(B)$ be spectra of $B$. If $A$ be normal, then $\dist\{\s(A),\s(A+B)\}\le |B|$
(see [Ka,p.291]).

 The normal operator $\cos (z-\ve \mV )$ has the eigenvalues
$\cos (z-\ve\n_j), j=1,..,N$ with the
multiplicity 2. Using the result from [Ka] and asymptotics
\er{312} and identity \er{318} we deduce that the eigenvalues $\D_j(z)$ of matrix $L(z)$ satisfy the asymptotics \er{319}. 
The proof of \er{aD2} is similar.
\BBox

\section {Proof of the main theorems}
\setcounter{equation}{0}

We need the following results from [K4].

\begin{lemma}  \lb{T41}
Let $V,V'\in \mH$. Then the following asymptotics hold true:
\[
\lb{31} \F(z,\n)=(\cos z-\n)^{2N}+o(e^{N|\Im z|}) \ \ \ \ \ as \ \ \ \ \ |z|\to \iy,
\]
where  $|\n|\le A_0$ for some constant $A_0>0$. Moreover, there exists an integer $n_0$ such that:

\no i) the function $\F(z,1)$ has
exactly $N(2n_0+1)$ roots, counted with multiplicity, in the disc
$\{|z|<\pi(2n_0+1)\}$ and for each $|n|>n_0$, exactly $2N$ roots,
counted with multiplicity, in the domain $\{|z-2\pi
n|<{\pi\/2}\}$. There are no other roots.

\no ii) the function $\F(z,-1)$
has exactly $2Nn_0$ roots, counted with multiplicity, in the disc
$\{|z|<2\pi n_0\}$ and for each $|n|>n_0$, exactly $2N$ roots, counted with multiplicity, in the domain $\{|z-\pi (2n+1)|<{\pi\/2}\}$. There are no other roots.

\no iii) Let in addition 
 $\n_i\ne \n_j$ for all $i\ne j\in \o_s$ for some $s=1,..,N_0$. Then the function $\r_s$
has exactly $2N_s(N_s-1)n_0$ roots, counted with multiplicity, in the disc $\{|z|<\pi (n_0+{1\/2})\}$ and for each $|n|>n_0$, exactly $N_s(N_s-1)$ roots, counted
with multiplicity, in the domain $\{|z-\pi n|<{\pi\/2}\}$.
There are no other roots. Moreover, 
\[
\lb{T32-2}
\r_s(z)=c_s\rt({\sin z+o(e^{|\Im z|})\/2z}\rt)^{N_s(N_s-1)},
\ \ \ \ \ c_s=\prod _{j,k\in \o_s}(\n_j-\n_k)^2, \ \ \  |z|\to\iy.
\]
\no ii) All zeros of $\r_s$ are given by $z_\a^{a,\pm}, \a=(j,k), j<k, j,k\in\o_s$ and $n\in \Z\sm \{0\}$. Furthermore, they satisfy
\[
\lb{T32-3}
z_\a^{n\pm}=\pi n+{\n_j+\n_k+o(1)\/2\pi n},\ \ 
\ \ \a=(j,k), \ \ n\to \pm\iy.
\]

\end{lemma}

\no {\bf Proof of Theorem 1.2.} i) 
We determine asymptotics \er{T2-1} for $z_j^{n\pm}$ as $n\to \pm\iy, j=1,..,N$. Lemma \ref{T32} yields $|z_j^{n\pm}-\pi n|<{\pi\/2}$, as $n\to \iy,j=1,2,..,2N$. Lemma \ref{T31} gives 
$
\D_m(z)=\cos (z-{\n_m\/2z})+O(1/z^2),m=1,..,N$ as   $z=\pi n+O(1)$.
For each $m=1,..,N$ there exists $j$ such that  $\D_j(z_m^{n\pm})=(-1)^n$. Thus we have $z_m^{n\pm}=\pi n+O(1/n)$.
Define the lokal parameter $\m$ by $z=\pi n+\ve\m, \ve={1\/2\pi n}$.
In order to improve these asymptotics
of $z_m^{n\pm}$ we need asymptotics of the $\P(1,z)$ as $z=\pi n+O(1/n), n\to \pm \iy$ given by \er{37}
$$
\P(1,z)=
e^{izJ_1}\rt(1-i\ve J_1\G_n+O(\ve^2)\rt),\ \ \G_n=\mV-iJ_1\wh {V_{n}'}, \ \ \ \ \ \wh {V_{n}'}=\int_0^1 V'(s)e^{i2\pi nJ_1s}ds
\ \ 
$$
where $\mV=\int_0^1 V^2(t)dt=\mV_0\os \mV_0$,\ \ 
$\mV_0={\rm diag} \{\n_1,...,\n_N\},\ \ \ \ 0\le \n_1\le \n_2\le...\le \n_N$. Thus we get 
\[
\lb{46} 
A={(-1)^n\P(1,z)-I\/-i\ve }={e^{i\ve \m J_1}\rt(I-i\ve J_1\G_n+O(\ve^2)\rt)-I\/-i\ve }
=J_1\rt(\G_n-\m+O(\ve)\rt).
\]
Hence  we study the zeros of the equation
\[
\lb{47} \det \rt(\mV-iJ_1\hat V_{(n)}'+O(\ve)-\m\rt)=0,\ \
\ \  \ \ \wh {V_{n}'}=\ma 0& \wh {v_n'}\\ \wh {v_n'}^*& 0\am,\ 
 \wh {v_n'}=\int_0^1 v'(t)e^{-i2\pi nt}dt, 
\]
where $ \m\in \C$.
We will use the standard arguments from the perturbation theory
(see [Ka,p.291]). Let $A,B$ be bounded operators, $A$ be a normal
operator and $\s(A),\s(B)$ be spectra of $A,B$. Then
$\dist\{\s(A),\s(B)\}\le\|A-B\|$.

Let $\z_m^{n\pm},(m,n)\in\{1,2,..,N\}\ts\Z$ be the eigenvalues of the
self-adjoint operator $\mV-iJ_1\wh{ V_{n}'}$.   Using the arguments from the perturbation theory (see [Ka,p.291]), we obtain that Eq. \er{47} has zeros $\o_m^{n\pm},(m,n)\in\{1,2,..,N\}\ts\Z$ such that
$\o_m^{n\pm}=\z_m^{n\pm}+O(n^{-1})$ as $n\to \pm\iy$,
which yields \er{T2-1}.

Consider the case $\n_{1}<...<\n_{N}$. We shall determine  asymptotics \er{T2-4} for the case $\a=(m,m), m=1,..,N$.
Let $z_\a^{n\pm}=\pi n+\ve\m$ and $\m-\n_m=\x\to 0$.
Using the simple transformation (unitary), i.e., changing
the lines and columns, we obtain
$$
\det \rt(\mV-iJ_1\wh {V_{n}'}+O(\ve)-\m\rt)=
\det \ma A_1& A_2\\ A_{3}& A_4\am=\det A_4\det K,
$$
$$
A_1=\ma -\x&b\\\ol b&-\x\am+O(\ve),\ \ \ \ \  b=-i \wh v_{n',\a}, \a=(m,m),
\  \ \ \ A_2,A_3=O(\d_n),\ \ \d_n=|\wh {V'_n}|+|\ve|,\  \ \
$$
$$ 
A_4=\mV_{0m} \os\mV_{0m}+O(\d_n),\ \ \ \ \ 
\mV_{0m}=\diag \{\n_j-\x, j\neq m\},
$$
$$
K=A_1-A_2A_4^{-1}A_3+O(\ve)=A_1+O(\f), 
\ \ \ \ \ \ \f=|\ve|+|\wh {V'_n}|^2,
$$
  which yields
$$
0=\det K=\x^2-|b|^2-\x b_1+ab_2+bb_3+O(\f^2), \ \ \ b_1,b_2,,b_3=O(\f),\ \ 
$$
where $b_1,b_2,b_3$ are analytic functions of $\x$. Rewriting
the last equation in the form $(\x+\a)^2=(c+\b)^2+O(\f^2), \a,\b=O(\f)$ and using the estimate $\sqrt{x^2+y^2}-x\le y$ for $x,y\ge 0$ we get 
$\x=\pm c+O(\f)$, which yields \er{T2-4} for $i=j$.

Consider the resonances. We shall determine  asymptotics \er{T2-4} 
for the case $\n_j\ne \n_{j'}$ for all $j\ne j'\in \o_s$ for some 
$s=1,..,N_0$. By Lemma \ref{T32}, the zeros of $\r_s$ have the form
$z_\a^{n\pm}, \a=(j,j'), j,j'\in \o_s, j<j', n\in \Z$
and satisfy $|z_\a^{n\pm}-\pi n|<\pi/2$. 

Asymptotics \er{aD2} yields
$\D_j(z)-\D_{j'}(z)=(\n_{j}-\n_{j'}){\sin z
\/2z}+O(z^{-2}e^{|\Im z|}),\ \ |z|\to\iy$. Then 
$   |z_\a^{n\pm}-\pi n|<\pi/2$ yields $z_\a^{n\pm}=\pi n+O(1/n)$ as $n\to \iy$.

We have
the identity   $\D_j(z)-\D_{j'}(z)=0$ at $z=z_\a^{n\pm}$.
Then using \er{319} we have
$$
\cos\lt(z_\a^{n\pm}-{\n_j\/2\pi n}\rt)- \cos\lt(z_\a^{n\pm}-{\n_{j'}\/2\pi n}\rt)=2(-1)^n\sin {\n_{j'}-\n_j\/4\pi n}\sin\lt(z_\a^{n\pm}-\pi
n-{\n_j+\n_{j'}\/4\pi n}\rt)={O(\d_n)\/n^2}
$$
which yields \er{T2-2}, i.e.,
\[
\lb{48}
z_\a^{n\pm}=\pi n+\ve(a_++\wt z_\a^{n\pm}),\ \
a_\pm={\n_j\pm \n_{j'}\/2},\ \ \ \ve={1\/2\pi n},\ \  \  \wt z_\a^{n\pm}=O(\d_n).
\]
We shall show that for large $n$ in the neighborhood of each
$\pi n+\ve a_+$ the function $(\D_j(z)-\D_{j'}(z))^2$
has two real zeros resonances (counted with multiplicity).
%Note that there are no others in the neighborhood of the point
%$\pi n$, since their number is equal to $N(N-1)$.
Introduce the functions
\[
\lb{49}
f_m(\m)=2(2\pi n)^2(1-(-1)^n\D_m(\pi n+\ve \m))=(\m-\n_m)^2+
O(\d_n).
\]
For the case $\m\to a_+$ we get
\[
\lb{410}
f_m(\m)=(a_+-\n_m)^2+o(1),  \ \ m=1,..,N, \ \ and \  \  
f_m(\m)=a_-^2+o(1),\ \ m=j,j'. 
\]
Hence the function $f_j-f_{j'}$ (maybe)
has the zeros, but the functions $f_j-f_m, m\ne j,j'$
have not zeros in the neighborhood of the point $a_+$.

Note that these functions are real outside the small 
neighborhood of $a_+$, otherwise for any complex branches
there exists a complex conjugate branch, but the asymptotics
\er{49} show that such branches are absent.

We have two cases: (1) let $f_s(\m),s=j,j'$
be real in some small neighborhood of $a_+$. 
Then the function $f_j-f_{j'}$ has at least
one real zero since by Theorem \ref{T1} the functions $f_j, f_{j'}$ are strongly monotone. Thus $(f_j-f_{j'})^2$ has at least 2 real zeros.

(2) Let $f_s(\m),s=j,j'$
be complex in some small neighborhood of $a_+$. Then they have
at least two real branch  points. Thus $(f_j-f_{j'})^2$ has at least 2 real zeros.

Hence $(f_j-f_{j'})^2$ has exactly two real zeros, since
the number of resonances (in the neighborhood of the point
$\pi n$) is equal to $N_s(N_s-1)$.

We determine the sharp asymptotics of resonances.
Recall that
$$
A={(-1)^n\P(1,z)-I\/-i\ve}=J_1(\mV-\m)-i\wh {V_{n}'}+O(\ve),\ \ 
\ \ \mV=\mV_0\os \mV_0.
$$
The operator $A-a_-$ has the eigenvalue $\x_0={(-1)^n\t_{n,s}-1\/i\ve}-a_-$
of multiplicity two, since $\t_{n,s}=(-1)^ne^{i\ve(a_-+o(1))}$.
The operator $J_1(\mV-a_+)-a_-=(\mV_0-\n_j)\os (\n_{j'}-\mV_0)$ has two eigenvalue ($= 0$) and other eigenvalues
are not zeros.
Using the simple transformation (unitary), i.e., changing
the lines and columns, $\m=a_++r\in \R, r=\wt r_{n,s}\to 0,$ we obtain
$$
\det (A-a_--\x )=
\det \ma A_1& A_2\\ A_3& A_4\am=\det A_4\det K,\ \
$$$$
 A_1=
\ma -r-\x & ib\\ i\ol b& r-\x\am +O(\ve),\ \ \  \  K=A_1-A_2A_4^{-1}A_3=\ma -r-\x+a_1& ib-\x+a_4\\ i\ol b+a_3& r+a_2\am
$$
$$
A_4=\mV_{0j}\os \mV_{0j'}+O(\d_n),
\ \ \ \ \ b=-\wh {v_{n,\a}'},\ \  \a=(j,j'),\ \ \   A_2, A_3=O(\d_n),
$$
the function $a_1,a_2,a_3,a_4=O(\f), \f=|\ve|+|\wh {V_{n}'}|^2$ and they analytic with respect to $\x$ in some small disk.
The function $\det K$ has the form 
\[
\det K=\x^2-r^2+|b|^2+a_1(r-\x)+a_2(-r-\x)-iva_3-ia_4\ol v-a_4a_3
\]
Then $0=\det K=(\x-\x_0)^2(1+O(\x-\x_0))$
 for $\x\to 0$ where $\x_0={((-1)^n\t_{n,s}-1)\/i\ve}-a_-$ is the zero of $F$ of multiplicity two.
Then $\x_0=O(\f)$ and we have $(r-\g)^2=(|b|-\b)^2+O(\f^2)$ where $\g,\b=O(\f)$. Then using the estimate $\sqrt{x^2+y^2}-x\le y$ for $x,y\ge 0$ we get $r=\pm |b|+O(\f)$.
\BBox

{\bf Proof of Theorem \ref{T3}} 
We consider $N_g^+$, the proof for $N_g^-$ is similar.

(i) Assume that $N_g^+=\iy$. Then, due to the Lyapunov-Poincar\`e Theorem and Theorem \ref{T2},
there exists a real sequence $z_k\to \iy$ as $k\to\iy$, such that
$z_k\in g^{n_k}_{\{j(m),m\}}$ for each $m=1,..,N$. Hence,
$\cap_{m=1}^Ng^{n_k}_{\{j(m),m\}}\ne \es$. 
Using  asymptotics \er{T2-4} and $k\to \iy$, we obtain $\n_1+\n_{j(1)}=...=\n_N+\n_{j(N)}$. Moreover, the estimates $\n_1<...<\n_N$ yield
$\n_{j(1)}^0>...>\n_{j(N)}^0$, i.e. $j(1)=N$, $j(2)=N-1$, ... Then,
$\n_1+\n_N=\n_2+\n_{N-1}=...$, which gives a contradiction.
\\
(ii) Let $2a=\n_1+\n_N=\n_2+\n_{N-1}=..$. Due to \er{T2-4}, $\pi n_k+a\in \cap_{m=1}^N g^{n_k}_{\{N+1-m,m\}}$ as $k\to\iy$. Then the Lyapunov-Poincar\`e Theorem yields $\pi n_k+a\notin\s(\cK)$, $k\to\iy$, i.e. $N_g=\iy$.
$\BBox$

\no {\bf Proof of Theorem \ref{T4}.} 
i) We need the following  result from [K4]:
 {\it Let $V\in \mH$. Define 
the quasimomentum $k_{j+N}=k_{j}=\arccos \D_j(z)=i\log \e(\D_j(z)), j=1,..,N$, see \er{dkm}, \er{dx}. 
 Then the averaged quasimomentum $k={1\/2N}\sum_1^{2N} k_j={1\/N}\sum_1^{N} k_j$ is analytic
in $\C_+$ and $k:\C_+\to k(\C_+)=\K(h)$ is a conformal mapping for
some $h\in C_{us}$. Furthermore, and there exist branches $k_j,j\in \ol{1,N}$ such that  \er{kc-1}-\er{T3-4} hold true}.

ii) Let $V, V'\in \mH$. 
We need the following results from [K4]: let
for some constants $C_0,C_1,C_2$ the following asymptotics hold
\[
\lb{410}
 \det (M(z)+M^{-1}(z))= \exp -{i2N\rt(z-{C_0\/z}-{C_1\/z^{2}}-
 {C_2+o(1)\/z^3}\rt)},\ \ \  {\rm  as}\ z=iy, \ \  y\to\iy.
\]
Then 
\[
 \lb{411}
k(z)=z-{Q_0\/z}-{Q_1\/z^2}-{Q_2+o(1)\/z^3},\ \ \  {\rm  as}\ y>r_0|x|, \ \  y\to\iy,\ \
\ \ {\rm for \ any} \ r_0>0,\ 
\]
where $ C_j=Q_j,\ j=0,1,2, \ \ Q_2=I_1^D+I_2^S-{Q_0^2\/2}$.
Using these results and  asymptotics from Lemma \ref{T42}  we obtain  \er{T4-1}-\er{T4-2}.\BBox

\section {Example of complex resonances}
\setcounter{equation}{0}

Let below $N=2$. Consider the operator $\cK_{\n,\t}=-iJ_1{d\/dt}+V_{\n,\t},\n=1, {1\/2}, {1\/3},.., \t\in \R$ acting in $L^2(\R)^4$, where the real periodic potential $V_{\n,\t}$ is given by
\[
\lb{Vtp} V_{\n,\t}=\ma \!\!0\!\!&\!\! v_{\n,\t}\!\!\\
\!\! v_{\n,\t}&\!\!0\!\!\am,\ \
v_{\n,\t}=-\ma a&\t b_\n(t)\\\t b_\n(t)&0\am,\qq {a\/2\pi}\in \R_+\sm \N,\qq
b_\n\in C(\T).
\]
 We need another representation
of $\cK_{\n,\t}$. Recall that $\mU={1\/\sqrt 2}(J_1+iJ)=\mU ^*$ and identities \er{222}
give $\cK_{\n,\t}^-=\mU \cK_{\n,\t}\mU=J{d\/dt}-V_{\n,\t}$.
%Recall that the operator $\cK y=Jy'+Vy, y=(y_1,y_2,y_3,y_4)^\top$.
Using the unitary transformation $y=(y_1,y_2,y_3,y_4)^\top\to \cU y=(y_1,y_3,y_2,y_4)^\top$ 
in $L^2(\R)^4$ we define the new operator $P_{\n,\t}=\cU \cK_{\n,\t}^-\cU^*={\bf j}{d\/dt}+W_{\n,\t}$, where 
\[
 W_{\n,\t}=-\cU V_{\n,\t} \cU^*=\ma a{\bf j}_2&\t b_\n{\bf j}_2\\ \t b_\n{\bf j}_2&0\am, \qq {\bf j_2}=\ma
0&1\\1&0\am,  \qq {\bf j}=\ma 0&1\\-1&0\am.
\]
We rewrite $W_{\n,\t}$ in the form $W_{\t,\n}=W^0+\t b_\n J_2$, where
$W^0=a\ma {\bf j_2}&0\\0&0\am$.
 If $\t=0$, then we have the unperturbed operator
$P^0=j{d\/dt}+W^0$ with a constant potential $W^0$.

{\bf 1. The $2\ts 2$ Dirac operator.} We consider the simple example of the $2\ts 2$ Dirac operator $P_1^0={\bf j}{d\/ dt}+a{\bf j_2}, a>0$ acting in $L^2(\R )\os L^2(\R )$.
The spectrum of $P_1^0$ is purely absolutely continuous and
consists of two intervals $(-\iy,-a)$,$(a,\iy)$ separated by the
gap $(-a,a)$. In this case only one gap is open and other gaps
are closed. The solution of the system
$jy'+a{\bf j_2}y=zy$  has the form $y=e^{\pm ikt}y_0,$ for
some constant vector $y_0\in \C^2 $ and $k$ satisfies
$$
\det (ik{\bf j}+a{\bf j_2}-zI_2)=-k^2-a^2+z^2, \ \ \ k=k(z)=\sqrt{z^2-a^2},
$$
where the quasimomentum $k:\C\sm [-a,a]\to \ \C\sm [ia,-ia]$ is
the conformal mapping with asymptotics $k(z)=z+o(1)$ as $|z|\to \iy$.
 Note that
\[
A={\bf j}(z-a{\bf j_2}),\ \ \ \
A^2=({\bf j}(z-a{\bf j_2}))^2=-z^2+a^2-za{\bf j}{\bf j_2}-za{\bf j}{\bf j_2}=a^2-z^2=-k^2.
\]
 Then the fundamental solution of the equation $j{\p_1^0}'+a{\bf j_2}\p_1^0=z\p_1^0, \p_1^0(0,z)=I_2 $ is given by
$$
\p_1^0(t,z)=e^{-tA}=\sum_{n\ge 0}{(-tA)^n\/n!}= \sum_{n\ge
0}\rt({(tA)^{2n}\/(2n)!}-{(tA)^{2n+1}\/(2n+1)!}\rt)
$$
\[
\lb{di1} =\sum_{n\ge0}\rt({(-1)^n(tk)^{2n}\/(2n)!}-
{A\/k}
{(-1)^n(tk)^{2n+1}\/(2n+1)!}\rt)=\cos tk-{A\/k}\sin tk.
\]
Thus the Lyapunov function has the form
\[
\lb{di2} \D_1^0(z)={\Tr e^{-A}\/2}=\cos k(z),\ \
k=k(z)=\sqrt{z^2-a^2},\ z\in\C\sm [-a,a].
\]

\no {\bf 2. The $4\ts4$ unperturbed operator.} Consider the $4\ts4$ operator $P^0={\bf j}{d\/ dt}+W^0$  in
$L^2(\R)^2\os L^2(\R)^2$, where  $W^0=\ma a{\bf j_2}&0\\0&0\am $
is the $4\ts4$ matrix.
We rewrite this operator in the form $P^0=P_1^0\os P_2^0$,  where the Dirac operators $P_1^0={\bf j}{d\/dt}+a{\bf j_2}$ and $P_2^0={\bf j}{d\/
dt}$ act in $L^2(\R)\os L^2(\R)$. The corresponding equation ${\bf j}{d\/
dt}\p^0+W^0\p^0=z\p^0$ has the fundamental solution $\p^0(t)$ given by
\[
\lb{up1} 
\p^0(t)=\p_1^0(t)\os\p_2^0(t),\ \ \ \ \
\p^0_1=e^{-At},\qq \p^0_2=e^{-{\bf j}zt},\ \ \ A={\bf j}(z-a{\bf j_2}),\ \ \ t\ge 0.
\]
Thus using \er{di1}-\er{up1}, \er{2l} we obtain that
$D^0(\t,\cdot)=\det (\p^0(1)-\t I_4)$ satisfies
$$
D^0(\t,\cdot)=
(\t^2-2\D_1^0\t+1)(\t^2-2\D_2^0\t+1),\qqq \D_m^0=(T_1^0-(-1)^m\sqrt{\r^0})/2,
$$
\[
\lb{up2} T_m^0(z)=\Tr (e^{-mA}+e^{-{\bf j}zm})=2(\cos
mk(z)+\cos mz), \qq  \r^0(z)=(\cos k(z)-\cos z)^2,
\]
for any $z\in \C, m=1,2$. We will determine the zeros of $\r^0(z)$. We
have $ 0=\cos k-\cos z=2\sin {k-z\/2}\sin {k+z\/2},\ \ \ \ k=k(z).
$ Then we obtain $k\pm z=2\pi n, n\in \Z\sm\{0\}$, which gives
zeros $r_{n}^0$ (each zero has multiplicity 2) of $\r^0$ by
\[
 \lb{up4}
 r_{n}^0=\pi n+{a^2\/4\pi n},\ \ \ \ k(r_{n}^0)=\cases
{\pi n-{a^2\/4\pi n}\ \ \ &if\ \ \ \ $|n|>{ a\/2\pi}$\cr
           -\pi n+{a^2\/4\pi n}\ \ \ &if\ \ \ \ $|n|<{ a\/2\pi}$\cr},
           \ \ \ n\in \Z\sm\{0\}.
\]
We determine the periodic spectrum for the equation ${\bf j}y'+W^0y=zy$.
Using \er{up1},\er{up2} we have
$
\det (\p^0(1,z)\mp I_4)=4(\cos k(z)\mp 1)(\cos z\mp 1),
$
which yields the periodic and anti-periodic spectrum multiplicity
2
\[
z_{n,1}^{0,\pm}=\pi n, n\in \Z, \ \ \ \ {\rm and} \ \ \
z_{n,2}^{0,\pm}=\pi n \sqrt{1+{a^2\/\pi^2 n^2}},\ \ n\in \Z\sm
\{0\},\ \ \
\]
and $z_{0,2}^{0,\pm}=\pm a$ has multiplicity one.
Note that the zeros $z_{n,p}^0\ne r_{m}^0$ for all $n,m\in
\Z,p=1,2$. In this case
there are no gaps in the spectrum. 

\no {\bf 3. The perturbed case.} 
Consider the $4\ts4$ operator $P_{\t,\n}={\bf j}{d\/ dt}+W_{\t,\n}$  in
$L^2(\R)^2\os L^2(\R)^2$, where   the $4\ts 4$ potential
$W_{\t,\n}=W^0+\t b_\n {\bf j_2} J_2$ satisfies \er{Vtp}, \er{Vn}.
We show that there exist the non-degenerated resonance gaps for
some $W_{\t,\n}$. In this example some resonances are real and
some are complex.
The fundamental solution $\p^{\t,\n}(t,z)$ of the Eq.
${\bf j} \p'+W_{\t,\n}\p=z\p$ satisfies the integral equation
\[
 \lb{pc1}
\p^{\t,\n}(t,z)=\p^0(t,z)+\t
\int_0^t\p^0(t-s,z){\bf j}b_\n(s){\bf j_2} J_2\p^{\t,\n}(s,z)ds.
\]
Then $\p^{\t,\n}(t,z)$ has asymptotics
\[
\lb{pc2} \p^{\t,\n}(m,z)=\p^0(m,z)+\t \p^1(m,z,\n)+\t^2
\p^2(m,z,\n)+O(\t^3e^{m|\Im z|})\qq as  \qq \t\to 0,
\]
uniformly in $\n=1,{1\/2},..,m=1,2, z\in \C$, where
\[
\lb{pc3} \p^1(m,z,\n)=\int_0^m\p^0(m-t,z)b_\n(t){\bf j_1}J_2\p_0(t,z)dt,
\]
\[
\lb{pc4}
\p^2(m,z,\n)=\int_0^m\p^0(m-t,z)b_\n(t){\bf j_1}J_2dt\int_0^t\p^0(t-s,z)
b_\n(s){\bf j_1}J_2\p_0(s,z)ds.
\]
The identity $\Tr \p^0(t,z){\bf j_1}J_2=0,\ \ (t,z)\in \R\ts \C$ yields
$\Tr\p^1(t,z,\n)=0$. Thus we obtain
\[
\lb{pc5} T_m^{\t,\n}(z)\ev \Tr\p^{\t,\n}(m,z)=T_m^0(z)+\t^2
T_{m2}(z,\n)+O(\t^3e^{m|\Im z|}),\qq m=1,2,
\]
\[
\lb{pc6} T_{m2}(z,\n)=\!\!\int_0^m\!\!\!\!
b_\n(t)dt\int_0^tb_\n(s)F_m(t,s,z)ds,\ \
F_m(t,s,z)=\Tr{\bf j_1}J_2\p^0(y,z){\bf j_1}J_2\p_0(\z,z),
\]
uniformly in $\n=1,{1\/2},.., |\t|<1, z\in \C$, where
$y=t-s,\z=m-y.$ 

\begin{lemma}  \lb{Tem}
Let $V_{\t,\n}$ satisfy \er{Vt}, \er{Vn}. Then the following asymptotics hold true
\[
 \lb{7f1}
T_{12}(z,\n)={T_1^0(z)\/2}+o(e^{|\Im z|}),\qqq
T_{22}(z,\n)=T_{2}^0(z)+4\f(z)+o(e^{2|\Im z|}),
\]
\[
\lb{7f2} T_1^{\t,\n}(z)=T_1^0(z)(1+{\t^2\/2})+o(\t^2e^{|\Im z|}),
\]
\[
\lb{7f3} T_2^{\t,\n}(z)=(1+\t^2)T_2^0(z)+\t^24\f(z)+o(\t^2e^{2|\Im z|}),
\]
\[
\lb{7f4}
\r^{\t,\n}(z)={T_2^{\t,\n}(z)+4\/2}-{T_1^{\t,\n}(z)^2\/4}=(1+\t^2)\r^0(z)+\t^22(\f(z)-1)+o(\t^2e^{2|\Im z|}),
\]
\[
\lb{7f5}
\det (\p^{\t,\n}(1,z)\mp  I_4)
=(1+\t^2)D^0(\pm 1,z)+\t^2
\rt(T_1^0-1-{\f(z)-1\/2}+o(e^{2|\Im z|})\rt),
\]
as $\n\to 0$, uniformly on $|\t|<1, z\in \C$, and where
$\f(z)=\cos k(z)\cos z+{z\/k(z)}\sin z\sin k(z)$.

\end{lemma}
 \no {\bf Proof.} Using $b_\n\to \d_{per}=\sum\d(t-{1\/2}-n)$
 and $F_j(\cdot,\cdot)\in C(\R^2),j=1,2$, we obtain
\[
\lb{4u}
\int_0^1dt\int_0^tb_\n(t)b_\n(s)F_1(t,s)ds
={1\/2}F_1\lt({1\/2},{1\/2}\rt)+o(e^{|\Im z|}),\
\
\]
\[
\lb{4v}
\int_0^2dt\int_0^tb_\n(t)b_\n(s)F_2(t,s)ds={1\/2}F_2\lt({1\/2},{1\/2}\rt)
+{1\/2}F_2\lt({3\/2},{3\/2}\rt)+F_2\lt({3\/2},{1\/2}\rt)+o(e^{2|\Im z|})
\]
as $\n\to 0$. Thus, if $m=1$, then 
\er{4u} gives
$$
T_{12}(z,\n)={1\/2}\Tr{\bf j_1}J_2\p^0(0,z){\bf j_1}J_2\p_0(1,z)+o(e^{|\Im z|})=
{T_1^0(z)\/2}+o(e^{|\Im z|}),
$$
If $m=2$, then $F_2({1\/2},{1\/2})=F_2({3\/2},{3\/2})=\Tr\p_0(2,z)$ and
the identity $ J_2\p^0=\ma 0&\p_2^0\\ \p_1^0&0\am $ yields
$$
F_2\lt({3\/2},{1\/2}\rt)=\Tr {\bf j_1}\ma 0 \ \ \ \ \p_2^0(y)\\
\p_1^0(y)\ \ \ \ 0\am {\bf j_1}\ma 0&\p_2^0(\z)\\ \p_1^0(\z)&0\am=2\Tr{\bf j_1}\p_2^0(1,z){\bf j_1}\p_1^0(1,z)
$$
$$
=2\Tr{\bf j_1}(\cos z-{\bf j}\sin z)  {\bf j_1}(\cos k-{A\/k}\sin k)
=2\Tr(\cos z+{\bf j}\sin z)(\cos k-{A\/k}\sin k)
$$$$
=4\cos z\cos k
-2\Tr{\bf j}{A\/2k}\sin z \sin k=4(\cos z\cos k+{z\/k}\sin z \sin k)=4\f.
$$
Then \er{4v} gives
$$
T_{22}(z,\n)=\Tr\p_0(2,z)+4\f(z)=T_2^0(z)+4\f(z)+o(e^{2|\Im z|})\qq
as\qq \n\to 0.
$$
 Substituting \er{7f1} into \er{pc5}
 we obtain \er{7f2}, \er{7f3}. 
The asymptotics \er{7f2}, \er{7f3} imply
$$
\r^{\t,\n}(z)={T_2^{\t,\n}(z)+4\/2}-{{T_1^{\t,\n}(z)}^2\/4}=(1+\t^2)\r^0(z)+\t^22(\f(z)-1)+o(\t^2e^{2|\Im z|}),
$$
and $D^{\t,\n}(\pm 1,\cdot)=\det (\p^{\t,\n}(1,\cdot)\mp  I_4)=
\!\!(T_1^{\t,\n}-1)^2-\r^{\t,\n}\!\!$ satisfies
$$
D^{\t,\n}(1,z)=
%\!\!(T_1^{\t,\n}-1)^2-\r^{\t,\n}\!\!
%$$$$
\rt((1+{\t^2\/2})T_1^0(z)-1\rt)^2\!\!-(1+\t^2)\r^0(z)-{\t^2\/2}(\f(z)-1)
+o(\t^2e^{|\Im z|})
$$
$$
=D^0(1,z)+\t^2\rt(T_1^0(z)(T_1^0(z)-1)-\r^0(z)-{(\f(z)-1)\/2}+o(1)\rt)
%=(1+\t^2)D^0(1,\cdot)+\t^2\rt(T_1^0-1-{\f(z)-1\/2}+o(1)\rt),
$$
as $\n\to 0$, which gives \er{7f5}. The proof for $D^{\t,\n}(-1,\cdot)$
is similar. \BBox

\no {\bf Proof of Proposition \ref{T5}}. We have the simple
asymptotics
\[
\lb{p71} k(z)=z-{a^2\/2z}+O(z^{-3}),\ \ \ \ \ \cos k(z)=\cos
z+{a^2\/2z}\sin z+O(z^{-2}e^{|\Im z|}),
\]
\[
\lb{p72} \r^0(z)=(\cos k(z)-\cos z)^2=\rt({a^2\sin
z\/2z}\rt)^2+O(z^{-3}e^{2|\Im z|})
\]
as $|z|\to \iy$.
We take $\vk>0$ such that the disks $B_n=\{|z-r_n^0|<\vk\},
n\in\Z$ are not overlapping. Define a constants
$\e_n=\min_{|z-r_n^0|=\vk}|\r^0(z)|>0$. Thus using \er{7f5} we get
\[
\lb{p73}
|\r^{\t,\n}(z)-\r^0(z)|=\t^2O(1)={\t^2\/\e_n}\r^0(z)O(1),\ \ \ \
|z-r_n^0|=\vk\ \ as \qq \t\to 0,
\]
for each $n$. We also obtain for $|z|=\pi n_0+1$ the following estimates
\[
\lb{p74} |\r^{\t,\n}(z)-\r^0(z)|=\t^2e^{2|\Im z|}O(1)=\t^2 |4\sin
z|^2O(1)=\t^2|z|^2|\r^0(z)|O(1),\ \ n_0\to \iy.
\]
Thus we take large $n_0$ and sufficiently small $\t$ such that $
|\r(z)-\r^0(z)|\le {1\/2}|\r^0(z)|$ on all contours $|z|=\pi
(n_0+1)$ and $|z-r_n^0|=\vk, |n|\le n_0$. Then by the Rouch\'e
theorem, $\r$ has as many roots, counted with multiplicity, as
$\r^0$ in the disks $\{|z|<\pi \pi n_0+1\}, \{|z-r_n^0|<\vk\}$.
Since $\r^0$ has exactly one double root at $r_n^0,n\ne 0$, and
since $n_0$ can be chosen arbitrarily large, we deduce that in each
disk $\{|z-r_n^0|<\vk\}, 1\le |n|\le n_0 $ there exist two zeros
$r_{n,\n}^\pm (\t)$ of $\r^{\t,\n}$ for sufficiently small
$\t,\n$.

Consider the zeros $r_{n,\n}^\pm (\t)$ of $\r^{\t,\n}$ in the disk
$\{|z-r_n^0|<\vk\}$ for fixed $n, 1\le n\le n_0$.
The proof for the case $n<0$ is similar. Recall
$\f(z)=\cos k(z)\cos z+{z\/k(z)}\sin z\sin k(z)$.
 Consider $\f(z)-1$ at the point
$r_n^0=\pi n+x_n, n\ge 1$, where $x_n={a^2\/4\pi n}$.  Using  \er{up4} we obtain
\[
\lb{p75}
 \f(r_{n}^0)-1=\cos^2 x_n-1+{r_{n}^0\/k(r_{n}^0)}\sin^2
x_n=\lt({r_{n}^0\/k(r_{n}^0)}-1\rt)\sin^2 x_n>0,\qq if \qq 1\le n<{a\/2\pi },
\]
\[
\lb{p77}
 \f(r_{n}^0)-1=\cos^2 x_n-1-{r_{n}^0\/k(r_{n}^0)}\sin^2
x_n=-\rt(1+{r_{n}^0\/k(r_{n}^0)}\rt)\sin^2 x_n<0,\qq if \qq n>{a\/2\pi }.
\]
We rewrite the function $\r^{\t,\n}$ in the disk
$\{|z-r_n^0|<\vk\}$ in the form
\[
\lb{p78}
R(z,\t)\ev{\r^{\t,\n}(z)\/(1+\t^2)^2}=(z-r_n^0)^2f(z)+\t^2\f_1(z,\t),\
\ \ f(z)={\r^0(z)\/(z-r_n^0)^2},\qq z\in B_n,
\]
for sufficiently small fixed $\n,\t$. The functions $R(z,\t), f(z), \f_1(z,\t)$ are analytic in $(z,\t)\in B_n\ts \{|\t|<\ve\}$ for some small
$\vk,\ve>0$ and satisfy
\[
\lb{p79} f(r_n^0)>0, \ \ \f_1(r_n^0,0)=2(\f(r_n^0)-1)+o(1)\qqq
as \qq \n\to 0.
\]
 Applying the Implicit Function Theorem to $R(z,\t)=0$ and using
\er{p75}-\er{p79} we obtain a unique solution $r_{n,\n}^\pm(\t)$
of the equation $\F(r_{n,\n}^\pm(\t),\t)=0,\t\in
(-\t_0,\t_0), r_{n,\n}^\pm(0)=r_n^0,$ for some $\t_0>0$ and here $r_{n,\n}^\pm(\t)$ is an analytic function in $\{|t|<\t_0\}$ and satisfies
\er{p1} with $R_n=-2(\f(r_n^0)-1)/f(r_n^0)$. 
The proof of the statement ii) is similar.
  \BBox

\section {Appendix}
\setcounter{equation}{0}

%Define  functionals

\begin{lemma} \lb{T42}
Let $V,V'\in \mH$ and let $r>0$. Then for $y\ge r|x|, y\to \iy$  following asymptotics hold: 
\[
\lb{412} \Tr L_1(z)=\|V\|^2\sin z,
\]
\[
\lb{413} \Tr L_2(z)=\rt(i\cH_1-{\Tr \mV^2\/2}+i{\|V'\|^2+o(1)\/2z}\rt)\cos z,
\]
\[
\lb{414} \Tr \mV L_2(z)=\cos z\Tr \rt(G_1+G-{\mV^3\/2}\rt)+o(e^{|\Im z}|),
\]
%\begin{multline}
\[
\lb{415} \Tr L_3(z)=i\cos z\Tr \rt(\int_0^1V_t^4dt+
G_1+G-{\mV^3\/6}\rt)+o(e^{|\Im z}|),\qqq \qqq \qqq \qqq \qqq \qqq
\]
$$
\qqq \qqq G_1=\mV \int_0^1J_1V_t'V_tdt,\ \ \  G=\mV \int_0^1J_1u_t'u_tdt,\ \ \ 
u_t=\int_0^tV_s^2ds,
$$
%\end{multline}
\[
\lb{416} \det L(z)=2^{-2N}\exp -i\rt(2Nz-{\cH_0\/2z}
 -{\cH_1\/(2z)^2}-{\cH_2+o(1)\/(2z)^3}\rt),
\]
where $\cH_0=\Tr\int_0^1V_t^2dt,\ \ \ 
\cH_1=\Tr\int_0^1-iJ_1V_t'V_tdt,\ \ \ \cH_2=\Tr\int_0^1\rt(V_t'^2+V_t^4\rt)dt.
$
\end{lemma}
\no {\bf Proof.} Identity \er{315} implies \er{412}. 
We will show \er{413}. Asymptotics \er{312} yields
\[
\lb{417}
L_2=-i\sin z\int_0^1J_1V_t'V_tdt+i{e_1J_1\/2}(\mV\wh {V'}+\wh {V'}\mV)+L_{21}+L_{22},
\]
$$
L_{21}=-{1\/2}\int_0^1\rt(e_{1}u'u+e_{-1}uu'\rt)dt,\ \ \ 
L_{22}={1\/2}\int_0^1\rt(e_{1}f'f+
e_{-1}ff'\rt)dt, \ \ f_t=\int_0^tV_s'e_{2s}ds
$$
as $|z|\to \iy, y\ge r|x|$.
We determine $\Tr L_2$. Using \er{29} we have $\Tr e_1J_1(\mV\wh {V'}+\wh {V'}\mV)=0$ and then
\[
\lb{418}
\Tr L_2=-i\sin z \cH_1+\Tr (L_{21}+L_{22}).
\]
Using \er{212} we get
\[
\lb{419}
\Tr L_{21}=-{\Tr \/2}\int_0^1\rt(e_{1}u'u+e_{-1}uu'\rt)dt=
-{\cos z\/2}\Tr \int_0^1\rt(u'u+uu'\rt)dt=-{\cos z\/2}\Tr \mV^2.
\]
 Due to  \er{212}, \er{21} we obtain
\[
\lb{420}
\Tr L_{22}={\Tr \/2}\int_0^1\rt(e_{1}f'f+
e_{-1}ff'\rt)dt=\Tr\int_0^1\int_0^te^{izJ_1(1-2t+2s)}V_t'V_s'dtds
\]
and using
\[
\lb{421} \int_0^1\!\!\int_s^te^{i2z(t-s)}f_tr_sdtds
=\int_0^1 {if_t'r_tdt\/2z}+{o(1)\/z},\ \ 
\int_0^1\!\!\int_s^te^{-i2z(t-s)}f_tr_sdtds
={o(e^{|\Im z|})\/z},
\]
for $f,h\in L^2(0,1)$  and $\Tr J_1(V')^2=0$, we have 
%\[\lb{422} 
$\Tr L_{22}={i\cos z\/2z}\rt(\|V'\|^2+o(1)\rt),
$
which give \er{413}. We will determine $\Tr \mV L_2$. Using \er{417} we have
\[
\lb{423} 
\Tr \mV L_2=\Tr \mV \rt(-i\sin z \int_0^1J_1V'Vdt+L_{21}  \rt)
+o(e^{|\Im z|}),
\]
since $\Tr \mV e_1J_1(\mV \wh {V'}+\wh {V'}\mV ) =0$ and $L_{22}(z)=o(e^{|\Im z|})$. Due to $e^{izJ_1}=\cos z +iJ_1\sin z$
and $\int_0^1(u'u+uu')dt=\mV^2, \Tr \mV^3=0$ we get
$$
\Tr \mV L_{21}=-{\cos z\/2} \Tr \mV^3-i{\sin z \/2}\Tr \mV J_1 
\int_0^1(u'u-uu')dt=-{\cos z\/2} \Tr \mV^3-i\sin z G
$$
and together with \er{423} we obtain \er{414}, since $\sin z=i\cos +O(e^{-|\Im z|})$.

We will determine $\Tr L_3$. Using \er{313} we rewrite $\Tr L_3$ in the form
$$
\Tr L_3=\Tr \vk_3=\Tr (KV^2+K^2V+KVK+K^3)\p_0=\sum_1^4 A_{k},\ \  %(Kf)(t)=\int_0^t\!\!\!e^{i(t-s)zJ_1}W(s)f(s)ds,
$$
\[
\lb{424}
A_1=\Tr KV^2\p_0,\ 
A_2=\Tr K^2V\p_0,\ A_3=\Tr KVK\p_0,\ \ \ A_4=\Tr K^3\p_0,
\]
and recall $ (Kf)(t)=\int_0^te_{t-s}W_sf(s)ds,\ \ \ 
 W=-iJ_1V^2-V' $.
Using \er{29},\er{211}, \er{212} and $e_1=\cos z+iJ_1\sin z$  we have
$$
A_1=\Tr\int_0^1e_{1}W_tV_t^2dt=\Tr\int_0^1e_{1}(-iJ_1V_t^2-V_t')V_t^2dt=\Tr\int_0^1-iJ_1e_1V_t^4dt=
$$$$
=\Tr\int_0^1-iJ_1(\cos z+iJ_1\sin z)V_t^4dt=\sin z\Tr \int_0^1V_t^4dt.
$$
The similar arguments give
$$
 A_2=\Tr\int_0^1\!\!\!\int_0^te_{1-t+s}W_te_{t-s}W_sV_sdtds
=\Tr\int_0^1\!\!\!\int_0^te_{1-t+s}(iJ_1V_t^2+V_t')e_{t-s}
(iJ_1V_s^2+V_s')V_sdtds
$$
\[
\lb{425}
=\Tr\int_0^1\!\!\!\int_0^tiJ_1\rt(e_1V_t^2V_s'V_s-e_{1-2t+2s}V_t'V_s^3\rt)dtds
\]
and 
$$
A_3=\Tr\int_0^1\!\!\!\int_0^t
e_{1-t+s}W_tV_te_{t-s}W_sdtds=\Tr\int_0^1\!\!\!\int_0^t
e_{1-t+s}(iJ_1V_t^2+V_t')V_te_{t-s}(iJ_1V_s^2+V_s')dtds
$$
\[
\lb{426}
=\Tr\int_0^1\!\!\!\int_0^tiJ_1\rt(e_1V_t'V_tV_s^2+e_{1-2t+2s}V_t^3V_s'\rt)dtds.
\]
Summing \er{425}, \er{426} we get $A_2+A_3=F_0+F_1$, where
\[
\lb{427}
F_0=\Tr\!\!\!\int_0^1\!\!\!\int_0^t\!\!\! iJ_1e_{1}
\rt(V_t^2V_s'V_s+V_t'V_tV_s^2\rt)dtds=\Tr iJ_1e_{1}
\mV\!\!\!\int_0^1\!\!\!V_t'V_tdt
=i\cos z\Tr J_1 \mV\!\!\!\int_0^1\!\!\! V_t'V_tdt,
\]
since $iJ_1e_1=iJ_1\cos z-\sin z$ and  $\Tr \mV\int_0^1V_t'V_tdt=0$.
We will show
\[
\lb{428}
F_1=\Tr\int_0^1\int_0^1 iJ_1e_{1-2t+2s}\rt(-V_t'V_s^3+  V_t^3V_s'\rt)dtds=o(e^{|\Im z|}).
\]
We use the standard arguments.
If $V,V',V''\in \mH$, then integration by parts gives \er{428}.
If $V,V'\in \mH$, then there exists $P_h,P_h',P_h''\in \mH$ such that $\|V-P_h\|+\|V'-P_h'\|=h$ for some small
$h\ge 0$. Then $F_1=
o(e^{|\Im z|})(1+O(h))$, which yields \er{428}, since
$h$ is arbitrary small.

The similar arguments give
$$
A_4=\Tr\int_0^1\int_0^t\int_0^s
e_{1-t+p}W_te_{t-s}W_se_{s-p}W_pdtdsdp=F_2+F_3,
$$
where
$$
F_2=\Tr\int_0^1\int_0^t\int_0^p(iJ_1)e_{1}V_t^2V_s^2V_p^2dtdsdp,\ 
$$$$
F_3=\Tr\int_0^1\int_0^t\int_0^p(iJ_1)\rt(e_{1+p-s}V_t^2V_s'V_p'+
e_{1-2t+2p}V_t'V_s^2V_p'+e_{1-2t+2s}V_t'V_s'V_p^2\rt)dtdsdp
$$
Using $u_s=\int_0^sV_t^2dt$ and $\int_s^1V_t^2dt=\mV-u_s$ we obtain
\[
F_2=\Tr(iJ_1e_{1})\!\!\int_0^1\!\!\!\int_s^1\!\!u_t'u_s'u_sdtds
=\Tr(iJ_1e_{1})\int_0^1\!\!(\mV-u_s)u_s'u_sds=i\cos z\Tr G-{\sin z\/6}\Tr \mV^3,
\]
since
$$\Tr(iJ_1e_{1})\int_0^1u_su_s'u_sds=
\Tr {iJ_1e_{1}\/3}\int_0^1(u_s'u_s^2+u_su_s'u_s+u_s^2u_s')ds=\Tr {iJ_1e_{1}\/3}\mV^3=-{\sin z\/3}\Tr \mV^3
$$
and
$$
\Tr(iJ_1e_{1})\mV\int_0^1u_s'u_sds=i\cos z\Tr J_1\mV\int_0^1u_s'u_sds-\sin z \Tr\mV\int_0^1u_s'u_sds
$$
$$
=i\cos z\Tr G-{\sin z\/2} \Tr\mV\int_0^1(u_s'u_s+u_su_s')ds=
i\cos z\Tr G-{\sin z\/2}\Tr \mV^3,
$$
where we used: $\Tr J_1\mV^3=0$
and $\Tr ABC=\Tr ACB$ for real self-adjoint matrix
and real representations \er{223}, \er{227} of $V$. 
Using standard arguments (see the proof of \er{428}) we have
$F_3=o(e^{|\Im z|})$. 
Summing $A_1,..,A_4$ we have \er{415}.

We will determine \er{416}. Asymptotics \er{312} yields
\[
\lb{430}
{L\/\cos z}=I_{2N}+S,\ \
S=i\ve \mV+{\ve^2L_2\/\cos z}+{\ve^3L_3\/\cos z}+O(\ve^4),\ \
S=O(\ve)
\]
as $\Im z\to \iy$.  In order to use the identity
$$
\det (I+S)=e^\F,\  \  \ \ 
\F=\Tr S-\Tr {S^2\/2}+\Tr {S^3\/3}+O(\ve^4),\ \
|S|=O(\ve),
$$
we need the traces of $S^m, m=1,2,3$. Due to \er{430}, we
get
%\[\lb{431}
$\Tr {S^3\/3}=-i\ve^3\Tr {\mV^3\/3}+O(\ve^{4})$.
Using \er{430},\er{412}-\er{415}  we get
\[
\lb{432}
-\Tr {S^2\/2}=\Tr \rt(\ve^2{\mV^2\/2}-i\ve^3{\mV L_2\/\cos z}+o(\ve^3)\rt)=\Tr \rt(\ve^2{\mV^2\/2}-i\ve^3\rt( G_1+G-{\mV^3\/2} +o(\ve^3)\rt)\rt),
\]
\[
\lb{433}
\Tr S=i\ve \cH_0+\ve^2\Tr \rt(i\cH_1-\ve^2{\mV^2\/2}\rt)+
i\ve^3\rt(\cH_2+
G_1+G-{\mV^3\/6}+o(1)\rt),
\]
and summing \er{431}-\er{433} we get 
$
\F=i\ve \cH_0+i\ve^2 \cH_1+i\ve^3 \cH_2+o(\ve^3)
$, which yields \er{416}.
\BBox

\no {\bf References}
\small

\no [BBK] Badanin, A.; Br\"uning, J.; Korotyaev, E. The Lyapunov function for Schr\"odinger operators with a periodic $2\times2$ matrix potential. J. Funct. Anal. 234 (2006), 106--126

\no [CK] Chelkak, D.; Korotyaev, E.
Spectral estimates for Schr\"odinger operators with periodic matrix potentials on the real line. Int. Math. Res. Not. 2006, Art. ID 60314, 41 pp.

\no [CG] Clark, S.; Gesztesy, F. Weyl-Titchmarsh $M$-function asymptotics, local uniqueness results, trace formulas, and Borg-type theorems for Dirac operators. Trans. Amer. Math. Soc. 354 (2002), no. 9, 3475--3534

\no [CHGL] Clark S.; Holden H.; Gesztesy, F.; Levitan, B.
Borg-type theorem for matrix-valued Schr\"odinger and Dirac
operators, J. Diff. Eqs. 167(2000), 181-210

\no [DS] Dunford, N. and Schwartz, J.: Linear Operators Part II:
Spectral Theory, Interscience, New York, 1988

\no [Fo]  Forster, O. Lectures on Riemann surfaces. 
Graduate Texts in Mathematics, 81. Springer-Verlag, New York, 1991

\no [GL] Gel'fand I.; Lidskii, V. On the structure of the regions
of stability of linear canonical    systems of differential equations with periodic coefficients.    (Russian)
   Uspehi Mat. Nauk (N.S.) 10 (1955), no. 1(63), 3--40

\no [Ge] Gel'fand, I.  Expansion in characteristic functions of an equation with  periodic coefficients. (Russian) Doklady Akad. Nauk SSSR (N.S.) 73, (1950), 1117--1120

\no [GKM] Gesztesy, F., Kiselev A.; Makarov, K. Uniqueness
results for matrix-valued Schrodinger, Jacobi, and Dirac-type
operators. Math. Nachr. 239/240 (2002), 103--145.

%\no [GK] Gohberg, I.; Krein, M. Theory and applications of Volterra %operators in Hilbert space. Translated from the Russian by A. Feinstein. %Translations of Mathematical Monographs, Vol. 24 American Mathematical %Society, Providence, R.I. 1970 x+430 pp.

\no [Ka] Kato, T. Perturbation theory for linear operators.
Springer-Verlag, Berlin, 1995

\no [KK1] Kargaev, P.; Korotyaev, E. Inverse Problem for the Hill
Operator, the Direct Approach.  Invent. Math., 129(1997), no. 3,
567-593

\no [KK2] Kargaev P.; Korotyaev, E.
Effective masses and conformal mappings. Comm. Math. Phys. 169 (1995), no. 3, 597--625

%\no [KK3] Kargaev P.; Korotyaev, E.
%Identities for the Dirichlet integral of subharmonic functions
%from the Cartright class,  Complex Var. Theory Appl. 50 (2005), no. 1, 35--50

\no [K1] Korotyaev, E.  Marchenko-Ostrovki mapping for periodic
Zakharov-Shabat systems,   J. Differential Equations, 175(2001), no. 2, 244--274

\no [K2] Korotyaev, E.  Inverse Problem and Estimates for Periodic
Zakharov-Shabat systems, J. Reine Angew. Math. 583(2005), 87-115

\no [K3] Korotyaev, E. Metric properties of conformal mappings on the complex plane with parallel slits, Internat. Math. Res.
Notices, 10(1996), 493--503

\no [K4] Korotyaev, E. Conformal spectral theory for the monodromy matrix, preprint 2006

\no [K5]
 Korotyaev, E. Inverse resonance scattering on the real line, Inverse Problems, 21(2005), 1-17

\no [Kr] Krein, M. The basic propositions of the theory of
$\lambda$-zones of stability of a canonical system of linear
differential equations with periodic coefficients. In memory of A.
A. Andronov, pp. 413--498. Izdat. Akad. Nauk SSSR, Moscow, 1955

%\no [MV] Maksudov, F.; Veliev, O.  Spectral analysis of differential %operators with periodic matrix coefficients. (Russian) Differentsial'nye %Uravneniya 25 (1989), no. 3, 400--409, 
%547; translation in Differential Equations 25 (1989), no. 3, 271--277

\no [MO] Marchenko V., Ostrovski I. A characterization of the
spectrum of the Hill operator. Math. USSR Sb. 26, 493-554 (1975)

\no [Mi] Misura, T. Properties of the spectra of periodic and anti-periodic boundary value problems generated by Dirac operators. I,II, Theor. Funktsii Funktsional. Anal. i Prilozhen, (Russian), 30 (1978), 90-101; 31 (1979), 102-109

%\no [Mi2] Misura T. Finite-zone Dirac operators.  Theor. Funktsii
%Funktsional. Anal. i Prilozhen, (Russian), 33 (1980), 107-11.

\no [Po]  Potapov, V. The multiplicative structure of $J$-contractive matrix functions. (Russian) Trudy Moskov. Mat. Ob\v s\v c. 4, (1955). 125--236
%, translated in  Amer. Math. Soc. Transl. (2) 15 1960 131--243

%\no [PT] P\"oshel,J., Trubowitz, E.: Inverse spectral theory. Pure
%and Applied Mathematics, 130. Academic Press, Inc., Boston, MA,
%1987. 192 pp.

\no [RS] M. Reed ; B. Simon. Methods of modern mathematical physics. IV. Analysis of operators. Academic Press, New York-London, 1978

\no [YS] Yakubovich, V., Starzhinskii, V.: Linear differential
equations with periodic coefficients. 1, 2.
   Halsted Press [John Wiley \& Sons] New York-Toronto,
   1975. Vol. 1, Vol. 2

\no [Z] Zworski M., : Distribution of poles for scattering
on the real line, J. Funct. Anal. 73, 277-296, 1987

\end{document}